\documentclass[11pt]{article}
\usepackage{mathrsfs}
\usepackage{amsfonts}
\usepackage{amsthm}
\usepackage{mathrsfs}
\usepackage{amsmath, amssymb}
\usepackage{shortvrb,psfrag}
\usepackage{enumerate}
\usepackage{color}
\usepackage{graphicx,subfigure}
\usepackage{dsfont}
\usepackage{latexsym}
\usepackage{longtable}

\makeatletter

\newcommand{\Rmnum}[1]{\expandafter\@slowromancap\romannumeral #1@}

\newcommand{\ur}{\mathrm{\textbf{r}}}
\newcommand{\um}{\mathrm{\textbf{m}}}
\newcommand{\rr}{\mathrm{\textbf{r}}}

\newcommand{\udeltayier}{\delta_{ij}}

\newcommand{\udeltasansi}{\delta_{kl}}
\newcommand{\udxyi}{\frac{\ud}{\ud x}}
\newcommand{\udxer}{\frac{\ud^2}{\ud x^2}}

\newtheorem{athm}{\bf \t}[section]
\newenvironment{thm} [1] {\def\t{#1}\begin{athm} \bf \rm} {\end{athm}}
\newcommand{\bthm}{\begin{thm}}\newcommand{\ethm}{\end{thm}}

\newtheorem{remark}{Remark}[section]

\newtheorem{proposition}{Proposition}[section]

\newcommand{\beq}{\begin{equation}}
\newcommand{\eeq}{\end{equation}}
\newcommand{\ben}{\begin{eqnarray}}
\newcommand{\een}{\end{eqnarray}}
\newcommand{\beno}{\begin{eqnarray*}}
\newcommand{\eeno}{\end{eqnarray*}}
\newcommand{\bali}{\begin{aligned}}
\newcommand{\eali}{\end{aligned}}

\numberwithin{equation}{section}

\newcommand{\ve}{\varepsilon}

\newcommand{\ud}{\mathrm{d}}

\newcommand{\xx}{\mathbf{x}}

\newcommand{\nn}{\mathbf{n}}

\newcommand{\mm}{\mathbf{m}}

\newcommand{\II}{\mathbf{I}}

\newcommand{\BS}{{\mathbb{S}^2}}
\newcommand{\BR}{{\mathbb{R}^3}}

\begin{document}

\title{From microscopic theory to macroscopic theory: a systematic study on static modeling for liquid crystals}

\author{Jiequn Han$^\dag$,\, Yi Luo$^{\dag *}$, \, Wei Wang$^{\ddag}$ and Pingwen Zhang$^{\dag *}$\\[2mm]
{\small $^\dag$ LMAM \& School of  Mathematical Sciences, Peking University, China}\\
{\small $ ^\ddag$ Beijing International Center for Mathematical Research, Peking University, China}\\
{\small $^*$ E-mail: royluo\_rv@gmail.com,\, pzhang@pku.edu.cn}\\
}

\date{\today}

\maketitle

\begin{abstract}
In this paper, we propose a systematic way of liquid crystal modeling to
build connection between microscopic theory and macroscopic theory.
 {A new Q-tensor theory based on Onsager's molecular theory which leads to liquid crystals with certain shape has been proposed}. Making uniaxial assumption, we can recover the
Oseen-Frank theory from the derived $Q$-tensor theory,  {and the Oseen-Frank model coefficients can be examined}. In addition, the smectic-A phase can also be characterized by
the derived macroscopic model.
\end{abstract}

\tableofcontents

\section{Introduction}

Liquid crystal (LC) phases are mesomorphic states between ordinary liquid and crystal. The constituent
LC molecules translate freely as in a liquid while exhibiting some long-range order above a
critical concentration (lyotropic) or below a critical temperature (thermotropic). The anisotropic properties make LC suitable for a wide range
of commercial applications. However, the inability to sufficiently control the degradation of orientational
order that LC display in the liquid state remains a great loss in potentially more important applications.
This fact highlights the need of establishing a simple and comprehensive mathematical model to capture
main characteristics of different LC phases and describe phase transition and defects.


 {The static LC models can be classified into three levels}: the molecular models, the tensor models and the vector models. The first
kind is microscopic theory, while the other two are macroscopic theories.
We shall begin by briefly reviewing these models.

The molecular models are based on the statistical theories of LC. In these models, the molecule has
a continuous distribution of orientations which corresponds to the actual physical situation.
However, the statistical mechanics of LC is so difficult that even for the simplest physical
models, exact solution  is very hard to work out.
Onsager \cite{Ons49} discussed the
statistics of a hard-rod system, and used a variational approximation to deal with the
non-linear integral equations. By making an additional approximation of the uniaxial
mean-field, Maier and Saupe \cite{MS59,MS60}  {suggested an analytic} thermodynamic
potential. A similar mean-field approximation to the Maier-Saupe theory was  {presented} by Doi
\cite{Doi81}. Most of the subsequent studies  {are based on} the
Maier-Saupe potential, such as the McMillan model \cite{McM71}, the Marrucci-Greco theory \cite{MG91} etc.
The molecular models are established on sound physical
theories, but they are not sensitive to macroscopic properties.
Moreover, the molecular
 { models posses high dimensional problems}.
%

The tensor model, also called Landau-de Gennes theory \cite{dGP95}, is a phenomenological theory
which ignores the detailed nature of the interactions and the molecular
structure. The free energy of these models is expressed as a functional of the
tensor order parameter $Q$. This order parameter is preferred as it is a good
measurement of macroscopic properties and it covers a wider class of LCs besides
simple nematics.   {A} variety of different extensions of the Landau-de Gennes
model have been proposed to study the sophisticated LC phases including
the cholesterics and the blue phases \cite{Bra75, HS88}. In spite of its success, the Landau-de Gennes model might
involve  {nonphysical solutions}. For this,  {Ball and Majumdar \cite{BM09}}
suggested that a modification to the entropy terms should be made to yield physically meaningful solutions.
Another problem  {with} the tensor theory is that it involves many phenomenological
coefficients which are difficult to  {decide using} experimental results.

The vector models, initiated by Oseen-Frank \cite{Ose33} and extended by
Ericksen \cite{Eri61, Eri91},   {are based on} continuum theory which disregards
the details of the structure on the molecular scale. It describes a weakly
distorted system in which, at any point, the changes in density of the
liquid induced by a long-range distortion are very small and the local
optical properties are still those of a uniaxial crystal. In terms of
a vector field, such distorted state may be described entirely.  The distortion energy
of the vector model can be interpreted into three parts indicating three
typical deformations: pure splay, pure twist and pure bend.
The elastic constants for these three parts, denoted by $K_1,K_2,K_3$,
play an important role in LC modeling. There are plenty of experimental methods,
for instance the Frederiks method and the transition method,
designed to measure these three constants for typical LC molecules under
certain circumstances \cite{DLRV69, MS60, Zve37}. And theoretical
and numerical investigations  of $K_1,K_2,K_3$ are also
abundant \cite{AF88, LM86, SS04}. However, the interpretations in terms of basic physical measurements
for these elastic constants remain unclear.

 {Despite the extensive literature on the static LC modeling}, little work
has been done to analyze the relations between different theories, especially between the macroscopic
theories and the microscopic theories. The above mentioned shortcomings and unsolved problems
in these theories also motivate us to build a unified framework to connect these models.
In  {this} framework, a simple and comprehensive macroscopic model should
be a simplification of the corresponding microscopic model applying the same molecular interaction potential.
Viewed in this light, the problem of determining macro coefficients is essentially a problem of representing
these by the original physical measurements, and it is no longer difficult to decide the number of
independent coefficients in the model. During reducing the complicated model into its simplifications,
information might be ignored or added, which is responsible for the occurrence of nonphysical solutions.

 {Following} this spirit, we propose a systematic way of LC modeling to build connection
among the three kinds of models: Onsager's molecular theory, Landau-de Gennes
Q-tensor theory and Oseen-Frank theory. Starting with Onsager's work, we
generalize it to the inhomogeneous system to characterize the distortion
of orientation by choosing suitable interaction potential in integral form. By applying local
Taylor expansion, we can write the energy in differential form, which is similar
to Marrucci-Greco's work \cite{MG91}. Next, by using the Bingham closure and truncating
at the low order moment, a $Q$-tensor model is obtained.
In this $Q$-tensor model, the physical constraint on $Q$ is  {automatically} satisfied. In addition,
the coefficients are determined by the molecular model, and their meaning can be
apparently interpreted.

Another important advantage of the new $Q$-tensor model is that
the well-known Oseen-Frank model can be recovered by restricting the
density to be  {a} constant and $Q$ to be uniaxial.
We can also calculate the values of the elastic coefficients in Oseen-Frank energy, and
examine the relation among them. Compared with former calculations of $K_1,K_2,K_3$, our expressions are more complete and precise.

Moreover, a model to characterize the simple smectic LC can also
be constructed by our method when introducing higher kernel function moment.
The form of whole free energy in Q-tensor model is similar to that
for nematic modeling and the layer thickness $d$ need not a priori
in our model. Numerical experiments show that the optimal solutions
are quite physical.

{This paper consists of two primary components of LC investigations. First
is the modeling where new nematic and smectic tensor models respecting the
physical mechanism are derived from the molecular statistical theories.
Further, these two models are consistent with each other. Second is the
study on the relationships among the existing three-level LC theories
which leads to a systematic way to compatibly model different phases
for different shaped LC molecules.}

  {This} paper is organized as follows. In section 2, we illustrate generally how to
derive macroscopic model, such as $Q$-tensor model and vector model
from the famous Onsager theory. We will apply this procedure for LC
with shape of rigid rod to model nematic phase, in section 3. A new $Q$-tensor
theory is derived, and the celebrated Oseen-Frank model is recovered there.
In section 4, we use the same way, but truncate at higher order of derivatives,
to model the smectic phase.  {Numerical} results ensure that the {smectic-A}
phase is captured.
We give several concluding remarks in section 5. Some detailed calculations
involved in the paper are provided in the appendix.

\section{A systematic way of static modeling of liquid crystals }

In this paper, we focus on the static modeling of liquid crystals. As we
would not like to take the boundary effect into account, we let
$\Omega\subseteq\BR$ be a periodic box. {In addition, we just consider molecules with axial symmetry.
Therefore the spacial information of one molecule can be specified by a position and a direction.} Use $\xx\in\Omega$ to denote the material point and
$f(\xx,\mm)$ to represent the number density for the number of
molecules whose orientation is parallel to $\mm$ at point $\xx$.
We start from Onsager's theory:
\begin{align}\label{eq:FreeEnergy-mole}
F[f]=k_BT\int_{\Omega}\int_{\BS}f(\xx,\mm)(\ln{f(\xx,\mm)}-1)
+\frac{1}{2k_BT}{\bar{U}(\xx,\mm)}f(\xx,\mm)\ud\mm\ud\xx,
\end{align}
where $k_B$ is the Boltzmann constant, $T$ is the absolute
temperature, and the mean-field interaction potential {$\bar{U}$} is defined by
\begin{align*}
{\bar{U}(\xx,\mm)}=k_BT\int_{\Omega}\int_{\BS}G(\xx,\mm; \xx',\mm')f(\xx',\mm')\ud\mm'\ud\xx'.
\end{align*}
Here $G(\xx,\mm;\xx',\mm')$ is the interaction kernel between
two molecules in the configurations $(\xx,\mm)$ and
$(\xx',\mm')$. In general, $G$ is translation invariant and hence it can be written in the form \beno
G(\xx,\mm;\xx',\mm')=G(\rr;\mm,\mm'), \eeno
where $\rr=\xx'-\xx$. The first part in (\ref{eq:FreeEnergy-mole}) represents the
entropy, while the second part describes the interaction energy between each
pair of two molecules in the system.

Firstly, we give the following two assumptions:
\begin{itemize}
\item[H1.] The LC state is very close to the equilibrium. Hence, we expect that the single particle distribution function is a satisfactory
but approximate basis to describe the macroscopic properties of the motion.

\item[H2.] The LC is composed of neutral particles surrounded by
force fields of short range compared with the average distance separating
the particles, that is, the LC is quite diluted.
\end{itemize}

H1 is the static modeling hypothesis. While, H2 is the critical hypothesis
which enables us to take the second virial expansion into account.
In the case when $c(\xx)=\int_\BS f(\xx,\mm)\ud\mm$ is small, the second virial expansion
is valid and the corresponding free energy approximation can be expressed in the form
\begin{align*}
&F[f]=F_0+{k_BT}\int  f(\xx,\mm){\ln f(\xx,\mm)}~\mathrm{d}\mm\mathrm{d}\xx  \nonumber \\
&\qquad\qquad+\frac{1}{2}\int f(\xx,\mm)G(\xx,\mm,\xx',\mm')f(\xx',\mm')~\mathrm{d}\mm'\mathrm{d}\xx'\mathrm{d}\mm\mathrm{d}\xx.
\end{align*}
Here the pairwise kernel function is defined as the classical expression for the second virial coefficient:
\begin{equation}
{G(\textbf{r},\textbf{m},\textbf{m}') = 1 - \exp (-U(\textbf{r},\textbf{m},\textbf{m}')/k_BT). \label{kernel}}
\end{equation}

{Here $U$ is an intermolecular potential. We assume $U$ consists of the average of interaction
of every pairs of basic particles. Hence $U$ is determined by the shape of molecules and the interaction potential $V$ between two basic particles, which can be taken as hard-core potential,
Lennard-Jones potential or other forms. Once $V$ is decided,} one
can then use the Monte-Carlo algorithm to numerically compute the molecular
model. As the computational cost for molecular model is too high, we must
look into the properties and the leading order moments of the kernel function.

It should be pointed out that for the higher density case, the second virial approximation
would not be sufficient. In such case it might be better to refer to the
Carnahan-Starling theory \cite{CS69} or other high density correction theory.
However, this  {is beyond the scope of this paper} and here we just use (\ref{kernel}) as our kernel function.

In the first paper \cite{XZ13} of this series , the relationship between the symmetry of molecule and the
properties of the kernel function $G$ has been discussed.
However, it is still difficult to derive the explicit expression of $G$.
For this reason, we turn to look into its moments by employing its
symmetric property. Since a non-local mean-field molecular
interaction potential is employed in our model, the
{orientational distribution function} $f(\textbf{x}',\textbf{m}')$ should be approximated
by its finite-order Taylor expansion series with respect to $\textbf{x}'$ at $\textbf{x}$:
\begin{align}
f(\textbf{x}',\textbf{m}')&=f(\xx+\rr,\mm')\nonumber\\ \label{f-expan}
&=f(\textbf{x},\textbf{m}') + \nabla f(\textbf{x},\textbf{m}')
\cdot\rr + \frac12 \nabla^2 f(\textbf{x},\textbf{m}'):\rr^T\rr + \cdots.
\end{align}
Then {$\bar{U}(\xx,\mm)$} can be formally written as
\begin{align*}
{\bar{U}(\xx,\mm)}=&k_BT\int_{\BS}\int_\Omega G(\rr; \mm, \mm')\Big\{f(\textbf{x},\textbf{m}')
 + \nabla f(\textbf{x},\textbf{m}') \cdot\rr \\
&\qquad\qquad+ \frac12 \nabla^2 f(\textbf{x},\textbf{m}'):\rr^T\rr + \cdots\Big\}\ud\rr\ud\mm'.
\end{align*}
For given kernel form of $G(\rr;\mm',\mm)$, we calculate the moments:
\begin{align*}
&M^{(0)}(\mm',\mm)=\int G(\rr,\mm',\mm)\ud\rr,\\
&M^{(1)}(\mm',\mm)=\int G(\rr,\mm',\mm)\rr\ud\rr,\\
&M^{(2)}(\mm',\mm)=\int G(\rr,\mm',\mm)\rr\rr^T\ud\rr,\\
&\qquad\cdots. 
\end{align*}
Then we get
\begin{align*}
{\bar{U}(\xx,\mm)}=&k_BT\int_{\BS}\Big\{f(\textbf{x},\textbf{m}')M^{(0)}(\mm,\mm')
+ M^{(1)}(\mm',\mm)\cdot\nabla f(\textbf{x},\textbf{m}')\\
&\qquad\qquad+ \frac12 M^{(2)}(\mm',\mm):\nabla^2 f(\textbf{x},\textbf{m}') + \cdots\Big\}\ud\mm'.
\end{align*}
The energy $F[f]$ becomes
\begin{align}\label{eq:FreeEnergy-mole-expan}
F[f]=&k_BT\int_{\Omega}\int_{\BS}\bigg\{f(\xx,\mm)(\ln{f(\xx,\mm)}-1)+\frac{1}{2}
\int_\BS M^{(0)}(\mm,\mm')f(\xx,\mm')f(\xx,\mm)\ud\mm'\bigg\}\ud\mm\ud\xx\nonumber\\
&\quad+\frac{k_BT}{2}\int_\Omega\int_\BS\int_\BS f(\xx,\mm)M^{(1)}(\mm,\mm')\cdot\nabla f(\xx,\mm')\ud\mm'\ud\mm\ud\xx\nonumber\\
&\quad+\frac{k_BT}{4}\int_\Omega\int_\BS\int_\BS f(\xx,\mm)M^{(2)}(\mm,\mm'):\nabla^2 f(\xx,\mm')\ud\mm'\ud\mm\ud\xx+\cdots.
\end{align}
The first line is independent of space variation of the probability distribution
function $f$.  We call it {\it bulk energy}, and denote by $F_{\text{bulk}}$.
The remainders, which depend on space variation of $f$, is called {\it elastic energy},
and denoted by $F_{\text{elastic}}$.

Now, we are going to express the energy by the spherical moments of $f$, namely,
\begin{align}
\int_\BS \underbrace{\mm\otimes\mm\otimes\cdots\otimes\mm}_{k \text{ times}} f(\xx,\mm)\ud\mm.
\end{align}
However, it is better to use the $k$-th order symmetric $traceless$ tensor
\begin{align}
Q_k[f]\triangleq\int_\BS\Xi_k(\mm)f(\mm)\ud\mm,
\end{align}
where $\Xi_k(\mm)$ is the $k$-th order symmetric $traceless$ tensor defined on the unit sphere, whose expression
for lower order takes the following form
\begin{align*}
&\Xi_1(\mm)=\mm;\\
&\Xi_2(\mm)=\mm\otimes\mm-\frac{1}{3}\II;\\
&\Xi_3(\mm)_{\alpha\beta\gamma}=m_{\alpha}m_{\beta}m_{\gamma}-\frac{1}{5}\Big(m_\alpha\delta_{\beta\gamma}
+m_\beta\delta_{\alpha\gamma}+m_\gamma\delta_{\alpha\beta}\Big);\\
&\Xi_4(\mm)_{\alpha\beta\gamma\mu}={m}_{\alpha}{m}_{\beta}{m}_{\gamma}{m}_{\mu}
-\frac{1}{7}\Big(m_{\alpha}m_{\beta}\delta_{\gamma\mu}+
m_{\gamma}m_{\mu}\delta_{\alpha\beta}+m_{\alpha}m_{\gamma}\delta_{\beta\mu}
+m_{\beta}m_{\mu}\delta_{\alpha\gamma}\nonumber\\
&\qquad\qquad\qquad+m_{\alpha}m_{\mu}\delta_{\beta\gamma}+
m_{\beta}m_{\gamma}\delta_{\alpha\mu}\Big)+\frac{1}{35}\Big(\delta_{\alpha\beta}\delta_{\gamma\mu}
+\delta_{\alpha\gamma}\delta_{\beta\mu}+\delta_{\alpha\mu}\delta_{\beta\gamma}\Big).
\end{align*}
One can see Appendix for precise definition for general $k$.
All components of $\Xi_k(\mm)$ are functions of linear combination of the $k$-th order spherical harmonics.
Moreover, for axial symmetric function $f(\mm)=f(\mm\cdot\nn)$, we have that
\begin{align}
Q_k[f]=S_k[f]\Xi_k(\nn),\text{ where } S_k[f]=\int_\BS P_k(\mm\cdot\nn)f(\mm)\ud\mm,
\end{align}
and $P_k(x)$ is the $k$-th Legendre's polynomial.

To derive tensor models from corresponding molecular models,
we need to use $Q_k(x)$ to express the total energy.
Since it is unrealistic to recover $f$ by finite number of moments,
we need to make closure approximation. It is very important to
choose a proper closure.  {To date}, variety of closure methods have been proposed.
For instance, we have the quadric closure (Doi closure), two Hinch-Leal
closures and the Bingham closure, etc.  {\cite{DE88, CL98, CLF95, FCL98, GMD00, IKKO03}}.
Here the Bingham closure is strongly suggested
for nematic phase and smectic phase, for the reason that Bingham closure
guarantees the existence of minimizers of the free energy functional and
provides with more accurate solutions. Additionally, it has
the good property of ensuring energy dissipation in kinetic models.
For given $f$ satisfying
\begin{align}
\int_\BS f(\xx,\mm)\ud\mm=c(\xx),\quad\int_\BS (\mm\mm-\frac13\II)f(\xx,\mm)\ud\mm=c(\xx)Q(\xx),
\end{align}
the Bingham closure is to use
\begin{align}\label{eq:fb}
f_{Q}(\xx,\mm)=c(\xx)\frac{\exp(B_{Q}(\xx):\mm\mm)}{\int_\BS\exp(B_Q(\xx):\mm\mm)\ud\mm}
\end{align}
to replace $f$, where $B_Q$ is a symmetric traceless matrix satisfying
\begin{align}
\int_\BS(\mm\mm-\frac13\II)\frac{\exp(B_Q:\mm\mm)}{\int_\BS\exp(B_Q:\mm\mm)\ud\mm}\ud\mm=Q.
\end{align}
It can be proved that $B_Q$ can be uniquely determined for given symmetric
traceless $Q$, if all the eigenvalues of  $Q$ belong to $(-\frac13, \frac23)$.
The bulk energy is then approximated by
\begin{align}
F_{\text{bulk}}=&k_BT\int_{\Omega}\int_{\BS}\bigg\{f_{Q_2}(\xx,\mm)(\ln{f_{Q_2}(\xx,\mm)}-1)\nonumber\\
&\qquad\qquad+\frac{1}{2}\int_\BS M^{(0)}(\mm,\mm')f_{Q_2}(\xx,\mm')f_{Q_2}(\xx,\mm)\ud\mm'\bigg\}\ud\mm\ud\xx.
\end{align}
Note that the above energy can be viewed as a functional of $c(\xx)$ and $Q_2(\xx)$.

Next, we consider the part of elastic energy.
To derive a convenient macroscopic model, we should only take finite terms
in (\ref{f-expan}) (or (\ref{eq:FreeEnergy-mole-expan})) into account.
If we want to model the nematic phase, it is natural to neglect the terms
whose order of derivatives are greater than two. If one would like to consider
the smectic phase, it seems enough to keep only the terms
whose order of derivatives are not greater than four.

Generally, we can truncate at $2m$-th order of derivatives. Now it is needed
to express the following terms in the energy using the tensors $Q_k$:
\begin{align}
\int_\BS\int_\BS f(\xx,\mm)M^{(l)}(\mm,\mm')
f(\xx,\mm')\ud\mm'\ud\mm', \text{ for }1\le l\le 2m.
\end{align}
For this, we have to separate the variables of $\mm$ and $\mm'$ in
$M^{(l)}(\mm,\mm')$.  {Generally} speaking, it can not be done precisely.
The reason is that $M^{(l)}(\mm,\mm')$ contains some terms like $|\mm\times\mm'|$.
Therefore, to deal with them,
we treat them as functions of $x=\mm\cdot\mm'$, and
use polynomial expansion, such as Taylor expansion or Legendre polynomial expansion,
to approximate them. In this way, we can express the energy by the tensors
$Q_k$. Also, the Bingham closure is used  to represent all $Q_k$ by $Q_2$:
\begin{align}
Q_k[f]=Q_k[f_{Q_2}],
\end{align}
where $f_{Q_2}$ is determined by (\ref{eq:fb}).
Together with the bulk energy part, we obtain a total energy
in $Q$-tensor form. Relevant introductions of the $Q$-tensor can be found in \cite{MN04}.

To derive the vector model, we only have to set $Q_2$ uniaxial:
\begin{align}
Q_2=S_2(\nn\nn-\frac{1}{3}\II).
\end{align}
Then an energy in the form of $\nn$ could be derived. If we regard the density $c(\xx)$
and the order parameter $S_2$ as constants,
then we can recover the well-known Oseen-Frank energy. The three important
coefficients can be directly expressed by the molecular parameters.


\subsection{Summary of the three-scale Schema for LC modeling}
In the above derivation, the key interrelationships among the molecular models, the tensor models
and the vector models can be tied up into a three-scaled
schema. We give some remarks in summary.

\begin{figure}[h]
   \centering
   \includegraphics[width=12cm]{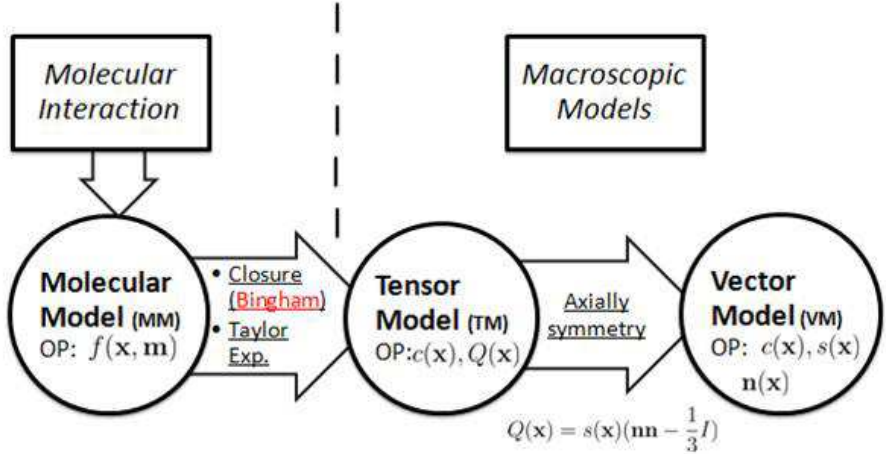}\\
   \caption{Three-scaled schema in static LC modeling.}
   \label{fig:schema}
\end{figure}

Firstly, {a pairwise kernel function}, which describes the intermolecular potential
and depends on concentration, temperature, molecular orientation and external field factors,
is the core element of a molecular model. The potential can generally be
classified into two categories: the lyotropic potential (such as the
hard-core potential) and the thermotropic potential (such as the
Lennard-Jones potential). When the intermolecular potential is properly
defined, the free energy functional of the molecular model can
be obtained immediately. The minimizers of the free energy functional are used to
describe the equilibrium state of the system.

Secondly, a molecular model can be changed into a tensor model through
Taylor expansion and closure approximation. The Taylor expansion is used
to collect moment information for the elastic energy part. And the closure
approximation will help to close the equations  {as} it will convert the
total free energy to a functional of the $Q$-tensor order parameter.
%

Finally, a tensor model can be easily changed into a vector model if
$Q$-tensor order parameter is restricted to uniaxial cases.
So far, we have elaborated the interrelationships of a three-scaled
schema, which can be further illustrated by Figure. \ref{fig:schema}.

The whole procedure can also be applied to molecules with complex shapes.
At that time, the selection of order parameters would be a basic and interesting problem.
The first paper \cite{XZ13} of this series discussed the relation between order parameters
and the molecular symmetry.

Furthermore, this three-scaled schema can also be employed in dynamic modeling
of LC system. One slight difference is that we need to additionally
guarantee the energy dissipation in characterizing dynamic fluids,
which might result in some difficulties in making closure approximation.
It can be proved that the Bingham closure  {\cite{Bin74}} satisfies the energy dissipation
law while the Doi closure  {\cite{DE88}} does not.
Another difference is that we can not derive the dynamical vector theory
(usually named Ericksen-Leslie theory) by simply setting $Q$  uniaxial.
Instead, it is needed to perform local expansion near the local equilibrium
as in the derivation from the Doi-Onsager theory to Ericksen-Leslie theory done in
\cite{KD83, EZ06,WZZ12}.
The third paper \cite{WZZ} in this series discussed how to employ this framework to
dynamical modeling for LC system.


\section{Modeling for nematic liquid crystals for rod-like molecule}

\begin{figure*}[h]
   \centering
   \includegraphics[width=5cm]{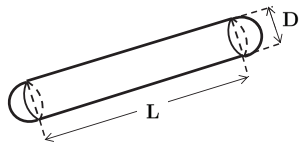}\\
   \caption{The geometry of the rod-like molecule}
   \label{fig:mol}
\end{figure*}

In this part, we will perform the above procedure to a particular
shape but widely studied molecule: the rod-like molecule.
The molecule is modeled as a round stick {with two caps}, see Fig. \ref{fig:mol}.

{This kind of molecule could be seen as a combination of spheres with same diameter $D$ alone a line with length $L$. Now we can write (\ref{kernel}) in the form
\begin{equation}
G(\textbf{r},\textbf{m},\textbf{m}') = 1 - \exp (-\frac{1}{k_BT L^2}
\int_{-L/2}^{L/2} \int_{-L/2}^{L/2} V(|\textbf{r} + t \textbf{m} -t'\textbf{m}'|) \,\mathrm{d}t \mathrm{d}t').
\end{equation}}To specify the kernel function as the second virial coefficient, we
have to decide potential $V(|\textbf{r}|)$ at first. The simplest
choice is the hard-core potential (or in the other name, the excluded-
volume potential), which was used in the Onsager theory. This potential
assumes molecules to possess hard elongated cores, which leads to the definition:
$$
 V = \left\{
    \begin{array}{ll}
    +\infty, & \text{if the two rods intersects};\\
    0, & \text{otherwise}.
    \end{array}
    \right.
$$
Apparently, the hard-core potential is a pure repulsive potential.
Another choice is the Lennard-Jones potential:
$$V^{(LJ)} (r) = 4 \frac{\ve}{k_B T} \{(\frac{\sigma}{r})^{12}- (\frac{\sigma}{r})^6\} $$
which takes the attractive interaction into consideration.

Notice that the hard-core potential is independent of temperature
$T$ while the Lennard-Jones potential possesses innegligible
 temperature
dependence. Because of this, the hard-core potential might be a proper
choice for modeling lyotropic LC while the Lennard-Jones fits better
for the thermotropic LC.  {As a matter of fact}, the actual temperature
for the nematic phase is often rather high, so both potentials give
qualitatively similar results in the nematic phase modeling.

For convenience, here we use the hard-core potential to demonstrate our further results. Such a steric repulsion gives
rise to steric cut-off effects. The kernel function $G$ takes the form:
\begin{equation}
G(\xx,\mm;\xx',\mm') = \left\{
    \begin{array}{ll}
    1, & \text{if the two rods intersects};\\
    0, & \text{otherwise.}
    \end{array}
    \right.
\end{equation}
It is needed to calculate $M^{(0)}(\mm',\mm),~
M^{(1)}(\mm',\mm),~M^{(2)}(\mm',\mm),$ et. al.
Of course, one can also consider the Lennard-Jones potential. But it would lead to
more complicate calculations.

Consider two rods with spherical ends pointing
$\textbf{m}$ and $\textbf{m}'$ respectively, the entire excluded
volume will be made up by three parts:
\begin{itemize}
\item region A (body-body): a $2D$-high parallelepiped whose section is a rhombus with side-length $L$ and angle $\gamma$;
\item region B (body-end): four semi-columns with side-length $L$ and radius $D$;
\item region C (end-end): four radius $D$ sphere at the corner.
\end{itemize}
\begin{figure*}[h]
   \centering
   \includegraphics[width=9cm]{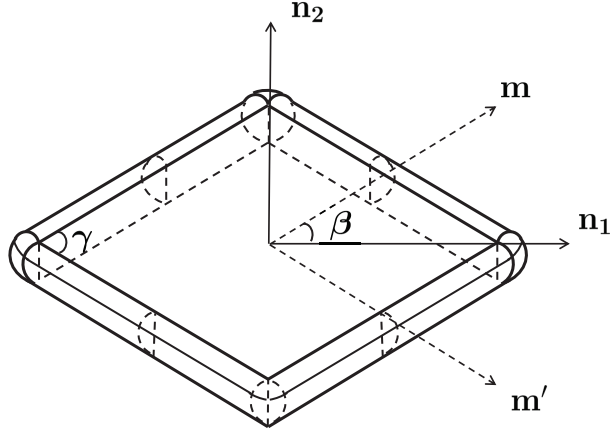}\\
   \caption{Hard-core potential: excluded-volume}
   \label{fig:lcn}
\end{figure*}
In the sequel, we introduce a dimensionless parameter $\eta (\le 1)$ as
\begin{align}
\eta=\frac{D}{L}.
\end{align}

\subsection{The bulk energy in $Q$-tensor form}
Calculating the volume in each of the above region, we can obtain the zero-th
moment of the kernel function (see Appendix for detailed calculation):
$$\int G(\textbf{r},\textbf{m},\textbf{m}')\mathrm{d}\textbf{r} =
2L^3\big(\eta\sin \gamma + \pi \eta^2 + \frac{2}{3}\pi \eta^3\big),$$
where $\sin \gamma = |\textbf{m}\times\textbf{m}'|$.
Thus, the bulk energy reads:
{\begin{align}\label{energy:bulk-onsager}
F_{\text{bulk}}[f] & = k_B T  \int_{\Omega} \int_{S^2}\left( f\ln f +  \int_{S^2}
L^3\big(\eta\sin \gamma + \pi \eta^2 + \frac{2}{3}\pi \eta^3\big)f(\xx,\mm')
\mathrm{d}\textbf{m}'\right)f(\xx,\mm)\mathrm{d}\textbf{m}\mathrm{d}\textbf{x}.
\end{align}}
Not surprisingly, this free energy functional is exactly the well-known
Onsager \cite{Ons49} model. Actually, Onsager integrated the hard-core potential
in a more complicated space field since he treated the LC molecule as
strict cylinders.

However, it is difficult to solve the minimizer problem of (\ref{energy:bulk-onsager}).
A possible way is to make a projection to the orthogonal polynomial space. In
this light, $|\textbf{m}\times\textbf{m}'|$ can be replaced by its
second order Legendre polynomial approximation
$-\frac{15\pi}{64}(\textbf{mm}-\frac13\II):(\textbf{m}'\textbf{m}'-\frac13\II)$.
Hence, the free energy functional becomes:
$$F_{\text{bulk}}[f]= k_B T  \int_{\Omega} \int_{\BS} f \mathrm{ln}
 f -\frac{15\pi L^3\eta}{64}\int_\BS (\textbf{mm}-\frac13\II):(\textbf{m}'\textbf{m}'-\frac13\II)ff'
 \,\mathrm{d}\textbf{m}'\mathrm{d}\textbf{m}\mathrm{d}\textbf{x},$$
which coincides largely with the Maier-Saupe \cite{MS59,MS60} model.
Actually, in the original presentation of Maier and Saupe, it
was assumed that intermolecular potential is due entirely to
van der Waals forces and is temperature-dependent. And it was
Doi's \cite{Doi81} work that used a mean-field approximation to  produce
the above functional and analyze bifurcations which occur as the
mean field is varied. In Doi's theory, $Q$-tensor defined as:
$$Q = \int_{\textbf{m}\in S^2} (\textbf{mm}-\frac13\II)\rho(\xx,\mm)\,\mathrm{d}\textbf{m}$$
is introduced (we always drop the subscript for brevity). Here,
\begin{align}
f(\xx,\mm)=c(\xx)\rho(\xx,\mm),\qquad c(\xx)=\int_\BS f(\xx,\mm)\ud\mm.
\end{align}
{To write the bulk energy in terms of $Q$, we use the Bingham closure as mentioned before.
For symmetric traceless matrix $Q$ whose eigenvalues belong to $(-\frac13, -\frac23)$,
let $B_Q$ be the unique symmetric traceless matrix such that
\begin{align}\label{eq:bq}
\int_\BS (\mm\mm-\frac13\II)\frac{\exp(B_Q:\mm\mm)}{\int_\BS\exp(B_Q:\mm'\mm')\ud\mm'}\ud\mm=Q.
\end{align}
Define
\begin{align}\label{eq:zq}
Z_Q=\int_\BS\exp(B_Q:\mm\mm)\ud\mm.
\end{align}
We replace the entropy term $\int_\Omega\int_\BS f\ln f\ud\mm\ud\xx$ by
\begin{align*}
\int_\Omega\int_\BS c(\xx)\frac1{Z_Q}\exp(B_Q:\mm\mm)\ln\Big(c(\xx)\frac1{Z_Q}\exp(B_Q:\mm\mm)\Big) \ud\mm\ud\xx,
\end{align*}
or equivalently,
\begin{align*}
\int_\Omega( c\ln c+c Q:B_Q-c\ln Z_Q)\ud\xx.
\end{align*}
Therefore, the energy can be simplified as follows:
\begin{align}\label{eq:FreeEnergy-mole-bulk}
F_{\text{bulk}} = k_B T  \int_{\Omega} c(\xx)\Big(\ln c(\xx)+Q(\xx):B_Q(\xx)-\ln
Z_Q(\xx) -\frac{15\pi L^3\eta}{64} c(\xx)|Q(\xx)|^2\Big)\mathrm{d}\textbf{x}.
\end{align}
}

The Maier-Saupe theory that we refer to here is temperature-independent.
We can also derive the Maier-Saupe theory from the Lennard-Jones
potential. The procedure is quite the same that if the kernel
function is based on the Lennard Jones potential, in the sense
of leading order, we have:
$$\int G(\mathbf{m},\mathbf{m'},\mathbf{r})~\mathrm{d}\mathbf{r}\approx H(L,D,T, \cos\gamma),$$
here $T$ is the temperature. Expanding $H(L,D,T, \cos\gamma)$ in
orthogonal polynomials with respect to the last variable, we can obtain
Maier-Saupe potential:
$$\int G(\mathbf{m},\mathbf{m'},\mathbf{r})~\mathrm{d}\mathbf{r}\approx H_1(L,D,T)-H_2(L,D,T)P_2(\cos{\gamma}).$$

\subsection{The elastic energy in $Q$-tensor form}

Now we turn to the elastic energy. Since $G(\rr; \mm',\mm)=G(-\rr; \mm',\mm)$, we have
$$M^{(1)}(\mm,\mm')=0.$$
To calculate the second moment, we introduce (for $\beta\neq0\text{~or~}\pi/2$)
$$
    \textbf{n}_1 = \frac{1}{2\cos \beta}(\textbf{m}+\textbf{m}'),\quad
    \textbf{n}_2 = \frac{1}{2\sin \beta}(\textbf{m}-\textbf{m}'),\quad
    \nn_3=\nn_1\times\nn_2.
$$
Under the coordinate $(\textbf{n}_1,\textbf{n}_2,\textbf{n}_3)$, the second moment must be a diagonal matrix:
$$M^{(2)} : = \int G(\textbf{r},\textbf{m},\textbf{m}')\textbf{r}\textbf{r}^T\,\mathrm{d}\textbf{r} = \mathrm{diag}(M_1,M_2,M_3).$$
Consequently, in the original coordinate, $M$ can be written as:
\begin{align*}
M^{(2)}  =&  M_1 \textbf{n}_1\textbf{n}_1 + M_2 \textbf{n}_2\textbf{n}_2 + M_3 \textbf{n}_3\textbf{n}_3\\
  = &M_3 \II + \frac{M_1-M_3}{4\cos^2 \beta}(\textbf{m}+\textbf{m}')(\textbf{m}+\textbf{m}') + \frac{M_2-M_3}{4\sin^2\beta}(\textbf{m}-\textbf{m}')(\textbf{m}-\textbf{m}')\\
  = & M_3 \II + (\frac{M_1}{4\cos^2\beta} + \frac{M_2}{4\sin^2\beta} - \frac{M_3}{\sin^2\gamma})(\textbf{mm}+\textbf{m}'\textbf{m}')\\
 &    + (\frac{M_1}{4\cos^2\beta} - \frac{M_2}{4\sin^2\beta} + \frac{M_3 \cos\gamma}{\sin^2\gamma})(\textbf{mm}'+\textbf{m}'\textbf{m})\\
  = & M_3 \II + (\frac{M_1}{4\cos^2\beta} + \frac{M_2}{4\sin^2\beta} - \frac{M_3}{\sin^2\gamma})(\textbf{mm}+\textbf{m}'\textbf{m}')\\
 &    + (\frac{M_1}{4\cos^2\beta \cos\gamma} - \frac{M_2}{4\sin^2\beta \cos\gamma} + \frac{M_3}{\sin^2\gamma})(\textbf{m}\cdot\textbf{m}')(\textbf{mm}'+\textbf{m}'\textbf{m})
\end{align*}
We can write:
\begin{align}\label{eq:M}
M^{(2)}=B_1\II + B_2(\textbf{mm}+\textbf{m}'\textbf{m}') + B_3(\textbf{m}\textbf{m}'+\textbf{m}'\textbf{m})(\mm\cdot\nn),
\end{align}
where $B_i$ are functions of $\gamma=\mm\cdot\mm'$:
$$
\left\{
\begin{array}{l}
B_1(\textbf{m}\cdot\textbf{m}') = M_3,\\
B_2(\textbf{m}\cdot\textbf{m}') = \frac{M_1}{4\cos^2\beta} + \frac{M_2}{4\sin^2\beta} - \frac{M_3}{\sin^2\gamma},\\
B_3(\textbf{m}\cdot\textbf{m}') = \frac{M_1}{4\cos^2\beta  \cos\gamma} - \frac{M_2}{4\sin^2\beta  \cos\gamma} + \frac{M_3}{\sin^2\gamma}.
\end{array}
\right.
$$
In the case of hard-core potential, it turns out that (see Appendix for details):
$$
\left\{
\begin{array}{ll}
B_1(\mm\cdot\mm')&=L^4D\Big(\frac{2|\mm\times\mm'|\eta^2}{3}
+\frac{\pi\eta^3}{2}+\frac{4\pi\eta^4}{15}\Big),\\
B_2(\mm\cdot\mm')&=L^4D\Big(\frac{|\mm\times\mm'|}{6}
+\frac{\pi\eta(1+\eta)}{3}+\frac{\pi\eta^3}{4}+\frac{2\eta^2}{3|\mm\times\mm'|}\Big),\\
B_3(\mm\cdot\mm')&=L^4D\eta^2\Big(\frac{2\arcsin(\mm\cdot\mm')}{3(\mm\cdot\mm')}
-\frac{2}{3|\mm\times\mm'|}\Big).
\end{array}
\right.
$$
It is worth pointing out  {here that} $B_1(\textbf{m}\cdot\textbf{m}')$, $B_2(\textbf{m}\cdot\textbf{m}')$,
$B_3(\textbf{m}\cdot\textbf{m}')$ are all even functions with respect
to $\cos \gamma = \textbf{m}\cdot\textbf{m}'$.

We want to use $c(\xx)$ and
\begin{align*}
Q_2(\xx)=\int \Xi_2(\mm)\rho(\xx,\mm)\ud\mm,\quad
Q_4(\xx)=\int \Xi_4(\mm)\rho(\xx,\mm)\ud\mm,
\end{align*}
to express the elastic energy (truncated to the second moment of the kernel function):
\begin{align}
F_{\text{elastic}}^{(2)}=\frac14\int_\Omega\int_\BS\int_\BS M^{(2)}:\nabla f(\xx,\mm')\nabla f(\xx,\mm)\ud\mm'\ud\mm\ud\xx,
\end{align}
where $M^{(2)}$ is defined by (\ref{eq:M}).

Let $x=\mm'\cdot\mm$, then we have the following Legendre polynomial expansion:
\begin{align*}
&|\mm\times\mm'|=\sqrt{1-x^2}=\frac{\pi}{4}-\frac{5\pi}{32}P_2(x)
-\frac{9\pi}{256}P_4(x)+\cdots,\\
&\frac{1}{|\mm\times\mm'|}=\frac{1}{\sqrt{1-x^2}}=
\frac{\pi}{2}+\frac{5\pi}{8}P_2(x)+\cdots,\\
&\frac{\arcsin{x}}{x}=\frac{\pi\ln2}{2}+\frac{5\pi}{16}(3-4\ln2)P_2(x)+\cdots.
\end{align*}
Hence, we can get
\begin{align*}
\frac{1}{L^4D}B_1(\mm\cdot\mm')&=\Big(\frac{\pi\eta^2}{6}+\frac{\pi\eta^3}{2}
+\frac{4\pi\eta^4}{15}\Big)-\frac{2\eta^2}{3}\cdot\frac{5\pi}{32}P_2(\mm\cdot\mm')
-\frac{2\eta^2}{3}\cdot\frac{9\pi}{256}P_4(\mm\cdot\mm')+\cdots,\\
\frac{1}{L^4D}B_2(\mm\cdot\mm')&=\frac{1}{6}\big(\frac{\pi}{4}
-\frac{5\pi}{32}P_2(\mm\cdot\mm')+\cdots\big)
+\frac{\pi\eta(1+\eta)}{3}+\frac{\pi\eta^3}{4}\nonumber\\
&\qquad+\frac{2\eta^2}{3}\big(\frac{\pi}{2}+\frac{5\pi}{8}P_2(\mm\cdot\mm')+\cdots)\nonumber\\
&=\frac{\pi}{24}+\frac{\pi\eta(1+2\eta)}{3}+\frac{\pi\eta^3}{4}
+(-\frac{5\pi}{3\cdot64}+\frac{5\pi\eta^2}{12})P_2(\mm\cdot\mm')+\cdots,\\
\frac{1}{L^4D}B_3(\mm\cdot\mm')&=\eta^2\Big(\frac{2}{3}
\big(\frac{\pi\ln2}{2}+\frac{5\pi}{16}(3-4\ln2)P_2(\mm\cdot\mm')+\dots\big)\nonumber\\
&\qquad-\frac{2}{3}\big(\frac{\pi}{2}+\frac{5\pi}{8}P_2(\mm\cdot\mm')+\cdots\big)\Big)\nonumber\\
&=\eta^2\Big((\frac{\pi\ln2}{3}-\frac{\pi}{3})
+(\frac{5\pi}{24}-\frac{5\pi\ln2}{6})P_2(\mm\cdot\mm')\Big)+\cdots.
\end{align*}
Denote
{\begin{align}
&\alpha_{11}=\frac{\eta^2}{6}+\frac{\eta^3}{2}
+\frac{4\eta^4}{15},\quad \alpha_{12}=-\frac{2\eta^2}{3}\cdot\frac{5}{32},\quad
\alpha_{13}=-\frac{2\eta^2}{3}\cdot\frac{9}{256},\nonumber\\
&\alpha_{21}=\frac{1}{24}+\frac{\eta(1+2\eta)}{3}+\frac{\eta^3}{4},\quad
\alpha_{22}=-\frac{5}{3\cdot64}+\frac{5\eta^2}{12},\nonumber\\
&\alpha_{31}=\eta^2\Big(\frac{\ln2}{3}-\frac{1}{3}\Big),\quad
\alpha_{32}=\eta^2\Big(\frac{5}{24}-\frac{5\ln2}{6}\Big).
\end{align}}
It is not hard to show that $P_n(\mm\cdot\mm')$ can be written as the tensor-inner product of $\Xi_{n}$, that is,
$$P_n(\textbf{m}\cdot\textbf{m}') = b_n\Xi_n(\textbf{m}):\Xi_n(\textbf{m}'),\quad \forall n\ge 1,$$
where $b_n$ is the highest order coefficient of $P_n(x)$: $b_1 = 1$, $b_2 = \frac32$, $b_3 = \frac52$,
$b_4 = \frac{35}{8},\cdots$.
Therefore, we can write
\begin{align*}
\frac{1}{\pi L^4D}B_1(\mm\cdot\mm')&=\alpha_{11}+ \frac32 \alpha_{12} \Xi_2(\mm):\Xi_2(\mm')
+\frac{35}{8}\Xi_4(\mm):\Xi_4(\mm')+\cdots,\\
\frac{1}{\pi L^4D}B_2(\mm\cdot\mm')&=\alpha_{21} + \frac32 \alpha_{22} \Xi_2(\mm):\Xi_2(\mm')+\cdots,\\
\frac{1}{\pi L^4D}B_3(\mm\cdot\mm')&=\alpha_{31} + \frac32 \alpha_{32} \Xi_2(\mm):\Xi_2(\mm')+\cdots.
\end{align*}
Dropping the high order terms in the above expansions, we arrive
\begin{align}
F_{\text{elastic}}^{(2)}
= & \frac{k_BT}{4} \int_\Omega\int_{\BS}\int_{\BS} M_{ij}^{(2)}\partial_i f(\textbf{x},\textbf{m})
\cdot\partial_j f(\textbf{x},\textbf{m}')
 \,\mathrm{d}\textbf{m}' \,\mathrm{d}\textbf{m}\ud\xx\nonumber\\
=&\frac12\int_\Omega\bigg\{J_1|\nabla{c}|^2+J_2|\nabla(cQ)|^2
+J_3|\nabla(cQ_4)|^2+J_4\partial_i(cQ_{ij})\partial_jc\nonumber\\
&+J_5\Big(\partial_i(cQ_{ik})\partial_j(cQ_{jk})+\partial_i(cQ_{jk})\partial_j(cQ_{ik})\Big)\nonumber\\
&+J_6\Big(\partial_i(cQ_{4iklm})\partial_j
(cQ_{4jklm})+\partial_i(cQ_{4jklm})\partial_j(cQ_{4iklm})\Big)\nonumber\\
&+J_7\partial_i(cQ_{4ijkl})\partial_j(cQ_{kl})\bigg\}\ud\xx.\label{Q-Energy}
\end{align}
The elastic coefficients can be written as  (see Appendix for details)
\begin{align*}
&J_1=-\frac{\pi}{2}{L^5\eta}k_BT\Big(\alpha_{11}+\frac23\alpha_{21}+\frac29\alpha_{31}+\frac{4}{45}
\alpha_{32}\Big),\\
&J_2=-\frac{\pi}{2}{L^5\eta}k_BT\Big(\frac32\alpha_{12}+\frac37\alpha_{22}
+\frac{9}{49}\alpha_{32}\Big),\qquad
J_3=-\frac{35\pi}{16}{L^5\eta}k_BT\alpha_{13},\qquad\\
&J_4=-\frac{\pi}{2}{L^5\eta}k_BT\Big(2\alpha_{21}+\frac{4}{3}
\alpha_{31}+\frac{2}{5}\alpha_{22}+\frac8{15}\alpha_{32}\Big),\\
&J_5=-\frac{\pi}{2}{L^5\eta}k_BT\Big(\alpha_{31}+\frac{25}{49}\alpha_{32}+\frac67\alpha_{22}\Big),\\
&J_6=-\frac{3\pi}{4}{L^5\eta}k_BT\alpha_{32},\qquad
J_7=-\frac{\pi}{2}{L^5\eta}k_BT\big(3\alpha_{22}+\frac{18}{7}\alpha_{32}\big).
\end{align*}
{It should be pointed out that since the hard-core molecular potential does not account for the temperature influence and we have
made several approximations and truncations in derivation, the expressions of the
coefficients $J_i$s may be not accurate.
What we want to suggest is the energy form (\ref{Q-Energy}).}
Now, we use the Bingham closure again to regard $Q_4$ as a tensor depending on $Q$:
\begin{align}
Q_4=\frac1{Z_Q}\int_\BS \Xi_4(\mm)\exp(B_Q:\mm\mm)\ud\mm,
\end{align}
where $B_Q$, $Z_Q$ are defined by (\ref{eq:bq}) and (\ref{eq:zq}). Then (\ref{Q-Energy}) is
a energy functional of $c$ and $Q$.
Together with the Maier-Saupe bulk energy part
\begin{align*}
F_{\text{bulk}}[c, Q] = k_B T  \int_{\Omega} c(\xx)\Big(\ln c(\xx)+Q(\xx):B_Q(\xx)
-\ln Z_Q(\xx) -\frac{15\pi L^3\eta}{64} c(\xx)|Q(\xx)|^2\Big)\mathrm{d}\textbf{x},
\end{align*}
we get the total free energy functional as:
\begin{align}\label{energy-QTensor}
F_\text{total}[c(\xx), Q(\xx)] = F_\text{bulk}[c(\xx), Q(\xx)] + F_\text{elastic}[c(\xx), Q(\xx)].
\end{align}
This is the main model hard-rods in our method.
We given some remarks in the following subsection.

\subsection{A brief look at our new $Q$-tensor model}

The energy (\ref{energy-QTensor})  {meets} the physical constraints
on the eigenvalues of the $Q$-tensor order parameter which guarantees
the existence of physically meaningful minimizers.{ All the model coefficients
are well interpreted in terms of the basic physical measurements and molecular
structures so that one can easily decide quantitatively proper values
for them in both numerical and physical experiments.} The concentration(or density),
is a spatially-dependent variable making great contributions to
modeling the LC smectic phase, which we will  {discuss later} in
this paper. In our framework, the $Q$-tensor models for the nematics
and the smectics are compatible with each other. They only differ in
the truncation process for the fourth-order moments terms.

For investigating the nematic phase only, the concentration variable
$c$ could be treated as a constant in the model. In this light, our
new $Q$-tensor model becomes the following form:
\begin{align}
F_\text{total}[Q(\xx)] = & F_\text{bulk}[Q(\xx)] + F_\text{elastic}[Q(\xx)]\nonumber\\
= & k_B T  \int_{\Omega} c\Big(\ln c+Q:B_Q
-\ln Z_Q -\frac{15\pi L^3\eta}{64} c|Q|^2\Big)\mathrm{d}\textbf{x}
\nonumber \\
&+ \frac12 c^2\int_\Omega\bigg\{J_2|\nabla Q|^2+J_3|\nabla Q_4|^2
+ J_5\Big(\partial_iQ_{ik}\partial_jQ_{jk}+\partial_iQ_{jk}
\partial_jQ_{ik}\Big)\nonumber\\
&+J_6\Big(\partial_iQ_{4iklm}\partial_j
Q_{4jklm}+\partial_iQ_{4jklm}\partial_jQ_{4iklm}\Big)\nonumber\\
&+J_7\partial_iQ_{4ijkl}\partial_jQ_{kl}\bigg\}\ud\xx
.\label{Q-Energy-1}
\end{align}

{This model is our modified version of the Landau-de Gennes tensor
model. Its bulk energy part will restrict the $Q$-tensor in the minimizers to
meet the physical constraints. And it also can be checked that if the values
of these measurements are picked in reasonable intervals, the elastic energy
will be bounded from below in any closed regions.
Another advantage of this new version is that no
phenomenological coefficients are involved in modeling. More detailed comparison
with the Landau-de Gennes model will be studied at the next subsection.}

If the axially-symmetry property is imposed on the minimizers of the total
free energy, which gives $Q(\xx) = S_2(\xx)(\textbf{n}(\xx)\textbf{n}(\xx)-\frac13 I)$,
the Ericksen vector model where
$$F_\text{total} = F_\text{total}[S_2(\xx), \textbf{n}(\xx)]$$
will be derived. Furthermore, if the scalar order parameter $S_2$ is regarded
spatially irrelevant, we finally arrive at the Oseen-Frank model where
 $$F_\text{total} = F_\text{total}[\textbf{n}(\xx)]$$
 The elastic coefficients $K_1, K_2, K_3$ can be interpreted then
 and their relations can be carefully examined.

 In addition to our derivation framework, similar analogue tensor models
 can be derived by following the same procedure when modeling for other
 molecular structures such as disk-shaped molecules and chiral molecules,
 or when the molecular interaction potential is redefined to account
 for the temperature dependence.

It is  {obvious} that after dropping some high order terms in (3.15),
we derive the Marrucci-Greco \cite{MG91} model which first analyzed the
long-range elasticity of LCPs. Some other nonhomogeneous extensions
of Doi's theory, differing in the intermolecular potential, can also
be obtained from our model by making truncations or approximations to
our nonlocal elastic potential. For instance, if we only keep the
$|\nabla Q_2|^2$ term in elastic energy, we get the Feng, Sgarlari
and Leal's \cite{FSL00} one-constant model; and if we slightly modify some
terms with $Q_4$, it turns out to be Yu-Zhang's \cite{YZ07} model which integrated
the long-range interactions in a ellipsoidal region. Besideds, the integral
form of our molecular model is quite the same as the model of Wang-E-Liu-Zhang
\cite{WELZ02}. Also, our method still works when the shape
of the molecules are changed, which will leads to Wang's \cite{Wan97} work.

\subsubsection{Comparisons with the Landau-de Gennes $Q$-tensor theory}
The well-known Landau-de Gennes model is a phenomenological theory, which has
successfully described the phase transition for liquid crystals.
Assuming that the free energy can be expanded as a power series of the order parameter $Q$ and of its spatial
derivatives, de Gennes gave the free energy functional as follows \cite{dGP95}:
\begin{eqnarray}\nonumber
F^{(LG)}[Q] &=& \int_{\Omega}\underbrace{\left(\frac{A(T-T^*)}{2} \mathrm{tr}(Q^2)
-\frac{B}{3}\mathrm{tr}(Q^3) + \frac{C}{4}(\mathrm{tr}Q^2)^2\right)}_{F^{(LG)}_{\text{bulk}}}\ud\xx\\
&&+\int_\Omega\underbrace{\left(L_1 Q_{ik,j}Q_{ij,k} + L_2 Q_{ij,j}Q_{ik,k} + L_3
|\nabla Q|^2 + L_4 Q_{lk}Q_{ij,k}Q_{ij,l}\right)}_{F^{(LG)}_{\text{elastic}}}\mathrm{d}\textbf{x}.\quad
\end{eqnarray}
The above energy contains two parts. The first part $F^{(LG)}_{\text{bulk}}$ governs the bulk effects,
$A,B,C$ are constants depending on temperature and material.
This expression of bulk energy is  {widely used} as it is capable of
describing a second-order phase transition and more importantly,
as it respects the axially-symmetry of the stationary points in homogeneous case.
The second part $F^{(LG)}_{\text{elastic}}$ is the elastic energy density that penalizes spatial non-homogeneities.
There are many works to study the solution of the classic
Landau-de Gennes model, for example,  one may see \cite{BM09, MZ09, Maj11} and the references therein.

An important problem in $Q$-tensor theory is to understand the physical
meaning of the tensorial order parameter $Q$. There is a kind of
interpretation that the $Q$-tensor represents the leading order
moment information of the {orientational distribution function}
$f(\textbf{x},\textbf{m})$. Hence, according to the derivation from the mean-field
approach, $Q$-tensor indicates the second moment tensor of $f$, i.e.
\begin{align*}
Q:=\int_{\BS}(\textbf{mm}-\frac13 \II)f~\mathrm{d}\textbf{m}.
\end{align*}
One can immediately draw a conclusion from the above definition that the molecular
theories require $Q$-tensor to be a symmetric, traceless $3\times 3$ matrix with
eigenvalues $\{\lambda_i(Q)\}$ constrained by the following inequalities:
$$-\frac13 \le \lambda_i(Q)\le \frac23, \qquad i=1,2,3.$$
Here, the inequalities were referred as physical constraints of the $Q$-tensor by
Ball and Majumdar \cite{BM09}, who pointed out that: {\it ``the bulk potential in the
Landau-de Gennes theory has no term that enforces the physical constraints on
the eigenvalues in the $Q$-tensor"}.

To avoid the non-physical flaws of Landau expansions,
the following bulk potential derived from the mean-field
Maier-Saupe free energy was suggested in \cite{BM09}:
\begin{align}F^{(BM)}_\text{bulk}=k_BT \inf\limits_{f \in \mathcal{A}_{Q}} \int_\BS
f(\textbf{m}) \mathrm{ln} f(\textbf{m}) \,\mathrm{d}\textbf{m} -\kappa |Q|^2,\label{energy:BM}
\end{align}
where
$$\mathcal{A}_Q = \Big\{f:S^2 \to \mathbb{R}, f \ge 0, Q=\int_\BS
(\textbf{mm}-\frac13\II)f(\textbf{m})\,\mathrm{d}\textbf{m},
\text{ and }\int_\BS f(\mm)\ud\mm=1. \Big\} .$$
They also proved the existence and the uniqueness of the solution
for the above functional.
It can be checked that
$$f_{B_Q}=\frac{\exp(B_Q:\mm\mm)}{\int_\BS\exp(B_Q:\mm\mm)\ud\mm}$$
solves the minimizing problem
$$\inf\limits_{f \in \mathcal{A}_{Q}} \int_\BS
f(\textbf{m}) \mathrm{ln} f(\textbf{m}).$$
Therefore, our bulk energy part (\ref{eq:FreeEnergy-mole-bulk})
is actually the same as the energy (\ref{energy:BM}).
As a result, the eigenvalues of the $Q$-tensor in our tensor model
are bounded both from below and above to meet physical constraints
due to their analysis in \cite{BM09}.

It is worth pointing out that by representing the $Q$-tensor as:
$$Q=s(\textbf{n}\textbf{n}-\frac13 I) + b(\textbf{n}'\textbf{n}'-\frac13 I), $$
the bulk energy density is in fact a function only of two scalars $s$ and $b$, i.e.
$$F^{(LG)}_\text{bulk}[Q]=\psi_\text{bulk}(s,b). $$
Indeed, we have that
$$\mathrm{tr}Q^2 =\frac23 (s^2 + b^2 -sb)  ,\quad
\mathrm{tr}Q^3 = \frac19 (2s^3 + 2b^3 -3s^2b -3sb^2) .  $$
Therefore, $F^{(LG)}_\text{bulk}[Q]$ is essentially the polynomials approximation
of the bulk energy density $F_{\text{bulk}}$ with
respect to $s$ and $b$. Therefore, the Bingham closure can provide us with an
approach to decide the material-dependent and temperature-dependent
coefficients $A,B,C$ in the Landau expansion of the bulk energy.

In addition, it has also been shown in \cite{BM09} that, for any boundary conditions,
if $L_4 \not= 0$, then the Landau-de Gennes energy $F^{(LG)}[Q] $
is unbounded from below. In other words, $\min_{Q} F^{(LG)}[Q] = - \infty$.
To modify the bulk energy part by (\ref{energy:BM}) or (regarding $c$ as a constant)
$$F_\text{bulk} = k_B T \bigg(B_Q:Q - \mathrm{ln}Z - \kappa |Q|^2\bigg),$$
leads a possible way to resolve this problem. However, it seems impossible
to prove the existence of global minimizers when $L_4 \not= 0$ unless suitable hypotheses
have been made on the elastic constants $L_i$. It might be widely
accepted that there are indeed some relationships on these elastic
constants. However, to the best of our knowledge, these relationships have
not been understood clearly yet.

By deriving the tensor model from the molecular scale, our work may offer
an approach to settle the above problem.  Since our tensor model is derived
from physical energy at molecular level, the elastic
energy should be bounded from below naturally.
On the other hand, when the density is assumed to be constant,
the elastic energy $F^{(2)}_\text{elastic}$ is a functional of $Q$.
Then by using suitable expansion, we can regard $F^{(LG)}_\text{elastic}$
as a certain approximation of $F_\text{elastic}(Q, \nabla Q)$,
and derive the elastic coefficients $L_i (1\le i\le 4)$ in terms of
the molecular parameters. This might give us a possible way  to
understand those coefficients.

\subsection{The elastic coefficients under uniaxial constrain}

On the vector scale, we will focus  on the famous Oseen-Frank
elastic constants $K_1, K_2, K_3$. Looking into their relations will
be a straight and important way to understand the nature of LC system.

\subsubsection{Derivation of the elastic coefficients}

In racemic or achiral system, the nematics show complete rotational symmetry
around their preferred orientation. Some previous investigators were then
inspired to restrict the LC model into uniaxial cases, which leads to the vector model.

One of the simplest and the most successful mathematical vector theory is
the Oseen-Frank theory \cite{Ose33} that characterizes the equilibrium by
a director field $\textbf{n}(\textbf{x}) \in \mathbb{S}^2$ with a spatially
invariant degree of orientational order. As a consequence of such assumptions,
the bulk potential in the Oseen-Frank model is a spatial-independent constant.
Therefore, the equilibrium configurations of LC are only relevant to the local
or global minimizers of the corresponding elastic energy functional.
The Oseen-Frank energy takes the form:
\begin{align}\nonumber
F_\text{elastc}^{(OF)}(\textbf{n}, \nabla \textbf{n}) = &\frac{1}{2}
K_1 (\mathrm{div} \textbf{n})^2 + \frac{1}{2}K_2 (\textbf{n}\cdot (\nabla
\times \textbf{n}))^2 + \frac{1}{2}K_3 |\textbf{n} \times (\nabla \times\textbf{n})|^2\\
& + \frac{K_2+K_4}{2} (\mathrm{tr}(\nabla \textbf{n})^2-(\mathrm{div} \textbf{n})^2),
\end{align}
where $K_i(i=1,2,3,4)$ are elastic constants, which describe three basic
types of distortions: pure splay, pure twist and pure bend respectively.
Apparently, the three elastic
constants $K_1, K_2, K_3$ should be non-negative, otherwise the energy
will be unbounded from below, and the global minimizer will be nonexistent.


To derive the celebrated Oseen-Frank vector model from tensor model, the following two constraints are required:
\begin{itemize}
\item $Q$-tensor should be restricted to the uniaxial form, i.e. $Q=S_2(\textbf{nn}-\frac13I)$;
\item The scalar order parameter $S_2$ should be spatially invariant.
\end{itemize}
It is worth pointing out that these two constraints are naturally satisfied in the equilibrium state of the homogeneous LC system.

We consider the case when the density $c$ is  constant.
Based on our previous analysis, the bulk energy functional can be written as
$$ F_\text{bulk}[Q] =  k_B Tc \int_{\Omega} \Big(B_Q:Q - \mathrm{ln}Z_Q - \frac{5\pi}{64}cL^2D |Q|^2 \Big)\mathrm{d}\textbf{x},$$
where $Z_Q=\int_\BS \exp(\mm\mm:B_Q)\ud\mm$.
One can easily figure out that the bulk potential of the Oseen-Frank
model actually only depends on the constant scalar order parameter $S_2$,
irrelevant of the orientation of molecules $\textbf{n}$. Particularly for the last term, we have
$$|Q|^2=\mathrm{tr}(Q^2)=\frac23 S_2^2.$$
Therefore, $F_\text{bulk}$ is an additive constant in the free energy
functional when $S_2$ is spatially invariant. Hence the bulk energy part can
be omitted in modeling. In other words, we need only concern about the
elastic potential of the system. The equilibrium state can be viewed
as the minimizers of the elastic energy functional in Oseen-Frank.

Writing the $Q$-tensor as
$Q=S_2(\textbf{nn}-\frac13 I)$ and utilizing the relation that $\partial_jn_in_i=0$,
the Oseen-Frank elastic constants can be derived as follows
(the detail calculation is included in Appendix):
\begin{align*}
&K_1=\pi c^2L^5\eta k_BT\Big(-S_2^2\Big(\eta^2\frac{299}{1568}-\frac{15}{7\cdot64}
-\frac{12\ln2\eta^2}{49}\Big)+S_4^2\Big(\frac{15\eta^2}{128}
-\frac{115\eta^2}{49}\Big(\frac{1}{8}-\frac{\ln2}{2}\Big)\Big)\nonumber\\
&\qquad\qquad+S_2S_4\frac{15}7\Big(-\frac{1}{64}+\frac{5\eta^2}{14}-\frac{3}{7}\eta^2{\ln2}\Big)\Big),\\
&K_2=\pi c^2{L^5\eta}k_BT\Big(-5S_2^2\Big(\eta^2\frac{19}{1568}-\frac{1}{7\cdot64}
-\frac{3\ln2\eta^2}{98}\Big)+S_4^2\Big(\frac{15\eta^2}{128}-\frac{15}{49}\eta^2\Big(\frac{1}{8}-\frac{\ln2}{2}\Big)\Big)\\
&\qquad\qquad+S_2S_4\frac{5}{7}\Big(-\frac{1}{64}+\frac{5\eta^2}{14}-\frac{3}{7}\eta^2{\ln2}\Big)\Big),\\
&K_3=\pi c^2{L^5\eta}k_BT\Big(-S_2^2\Big(\eta^2\frac{299}{1568}-\frac{15}{7\cdot64}
-\frac{12\ln2\eta^2}{49}\Big)+S_4^2\Big(\frac{15\eta^2}{128}
-\frac{150}{49}\eta^2\Big(\frac{1}{8}-\frac{\ln2}{2}\Big)\Big)\\
&\qquad\qquad-S_2S_4\frac{20}7\Big(-\frac{1}{64}+\frac{5\eta^2}{14}-\frac{3{\ln2}\eta^2}{7}\Big)\Big),
\end{align*}
where $S_4$ is a function of $S_2$ due to the Bingham closure, in other words, given by
\begin{align}
S_4=\int_\BS P_4(\mm\cdot\nn)\frac{\mathrm{e}^{r(\mm\cdot\nn)^2}}{\int_\BS \mathrm{e}^{r(\mm'\cdot\nn)^2}\ud\mm'}\ud\mm,
\end{align}
where $r=r(S_2)$ is the unique real number satisfying
\begin{align*}
S_2=\int_\BS P_2(\mm\cdot\nn)\frac{\mathrm{e}^{r(\mm\cdot\nn)^2}}{\int_\BS \mathrm{e}^{r(\mm'\cdot\nn)^2}\ud\mm'}\ud\mm,
\end{align*}
or equivalently
\begin{align}\label{eq:s2-1}
S_2=\frac{1}{2}\frac{\int_0^1(3z^2-1)\mathrm{e}^{r z^2}\ud z}{\int_0^1\mathrm{e}^{r z^2}\ud z}.
\end{align}

\begin{remark}
It can be proven that $S_2\in(-\frac{1}{2}, 1)$ is a monotonically increasing function of $r\in(-\infty, +\infty)$. Thus, $r$
can be uniquely determined for all $S_2\in(-\frac{1}{2}, 1)$.
\end{remark}

Next, based on the above expressions, we will go deep into the comparative relationships among $K_1, K_2, K_3$,
which are closely linked with the stability of the nematic phase.

\subsubsection{The comparative relationship among the Oseen-Frank elastic constants}

Although we have already expressed $K_1, K_2, K_3$ in terms of
$S_2, c, L, D, T$, it is still hard to decide their relationships.
Notice that there are only two independent variables which are:
$$\alpha = \pi L^2D c,\quad \mathrm{and} \quad \eta = D/L.$$
The problem lies in the fact that the scalar order parameter $S_2$ is in fact
related to other parameters. Indeed, the relation between $S_2$
and $\alpha$ is given by (\ref{eq:s2-1}) and (\ref{eq:s2-2}).
We will  {offer a simple explanation to this relation here}.
The order of coefficients of $F_\textrm{elastic}^{(2)}$ is $O(L^5)$, while
the order of $F_\textrm{bulk}$ is $O(L^3)$.
Hence, $S_2$ should be chosen to minimize $F_\textrm{bulk}$. By the work
of \cite{LZZ05} or \cite{FS05}, we know that
$S_2$ should be chosen  {to} satisfy
\begin{align}\label{eq:s2-2}
\frac{15\alpha}{32}\int_0^1(3z^2-1)\mathrm{e}^{r z^2}\ud z=2r{\int_0^1\mathrm{e}^{r z^2}\ud z}.
\end{align}
where $r=r(S_2)$.

Once the concentration parameter $\alpha$ and the diameter-to-length
ratio $\eta$ are decided, we can compute $K_1, K_2, K_3$ for lyotropic
LC. If we assume that the LC molecule has a diameter of $\sim 5$ \r{A} with
temperature $T=400 K$ and dimensionless volume fraction
$\Phi=c\cdot\frac14 \pi LD^2 = 0.4$, the elastic constants $K_i$s
will have the dimension of energy/cm (or dynes) and the expected
magnitude is $10^{-6}$ dyn. This is indeed the correct order of
magnitude according to the famous Frederiks experiment \cite{FZ33}
conducted by Prost et. al \cite{Zve37}, Saupe \cite{MS60} and Durand
et. al. \cite{DLRV69}.

\begin{figure}
   \centering
      \includegraphics[width=7cm]{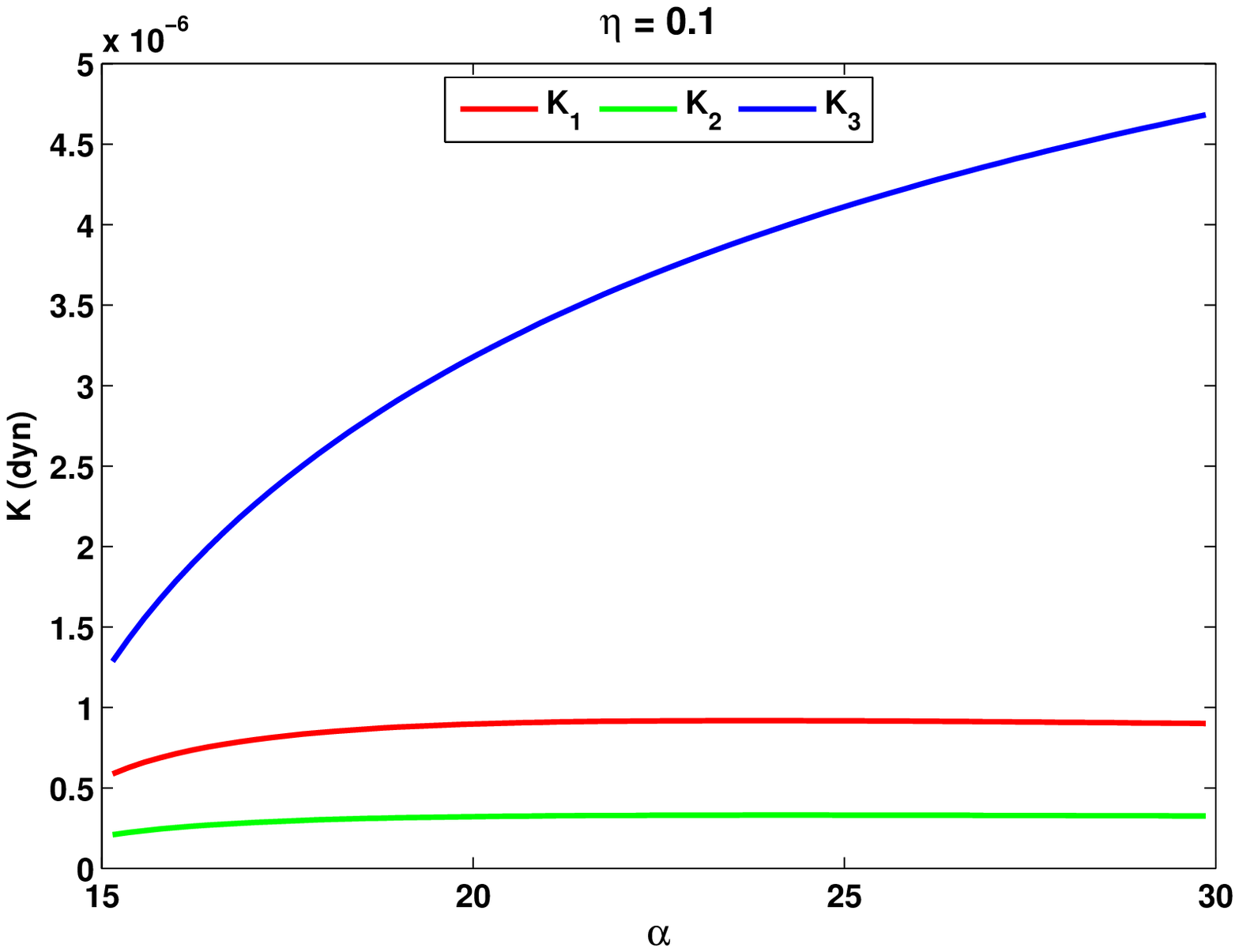}
      \includegraphics[width=7cm]{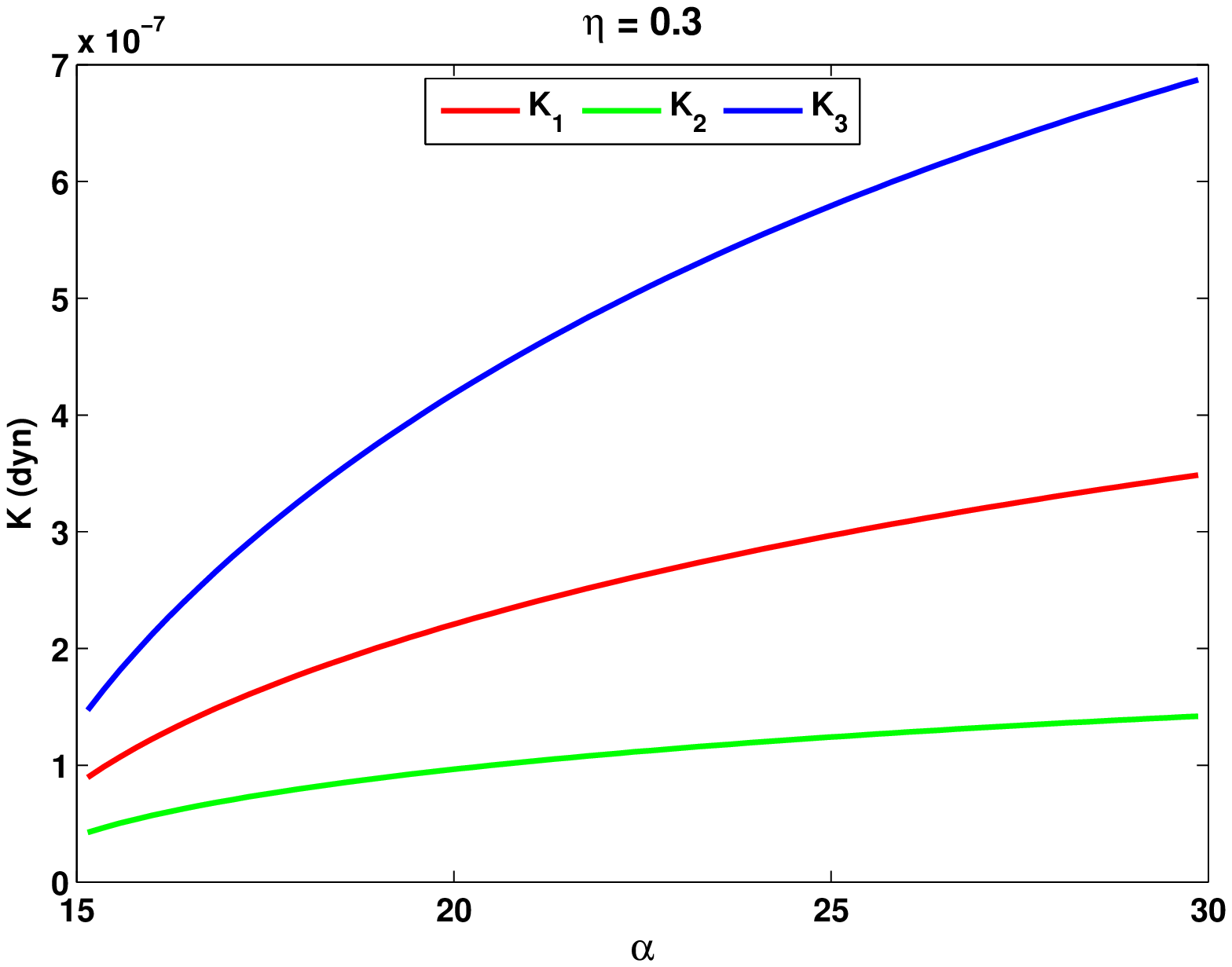}\\
         \centering
      \includegraphics[width=7cm]{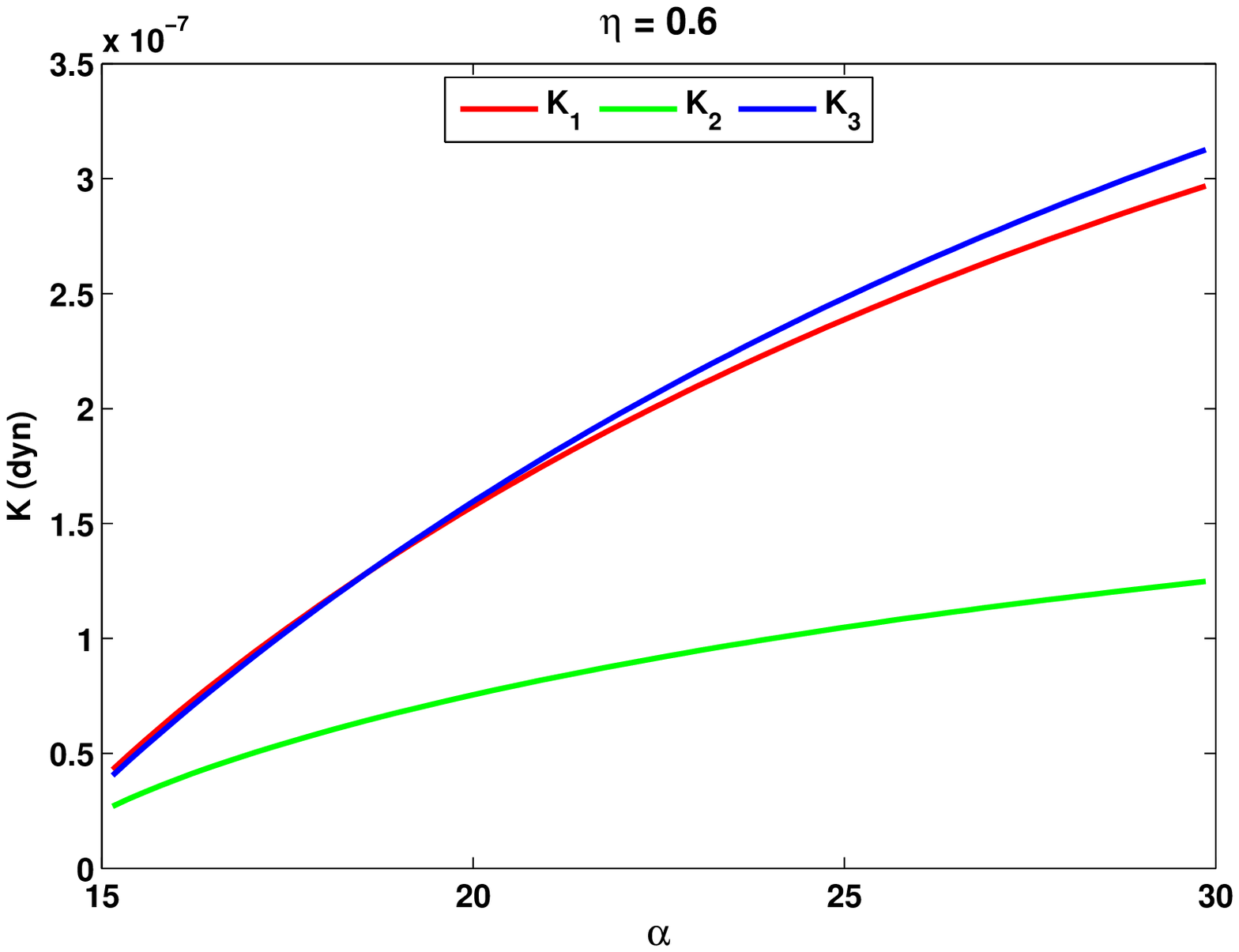}
      \includegraphics[width=7cm]{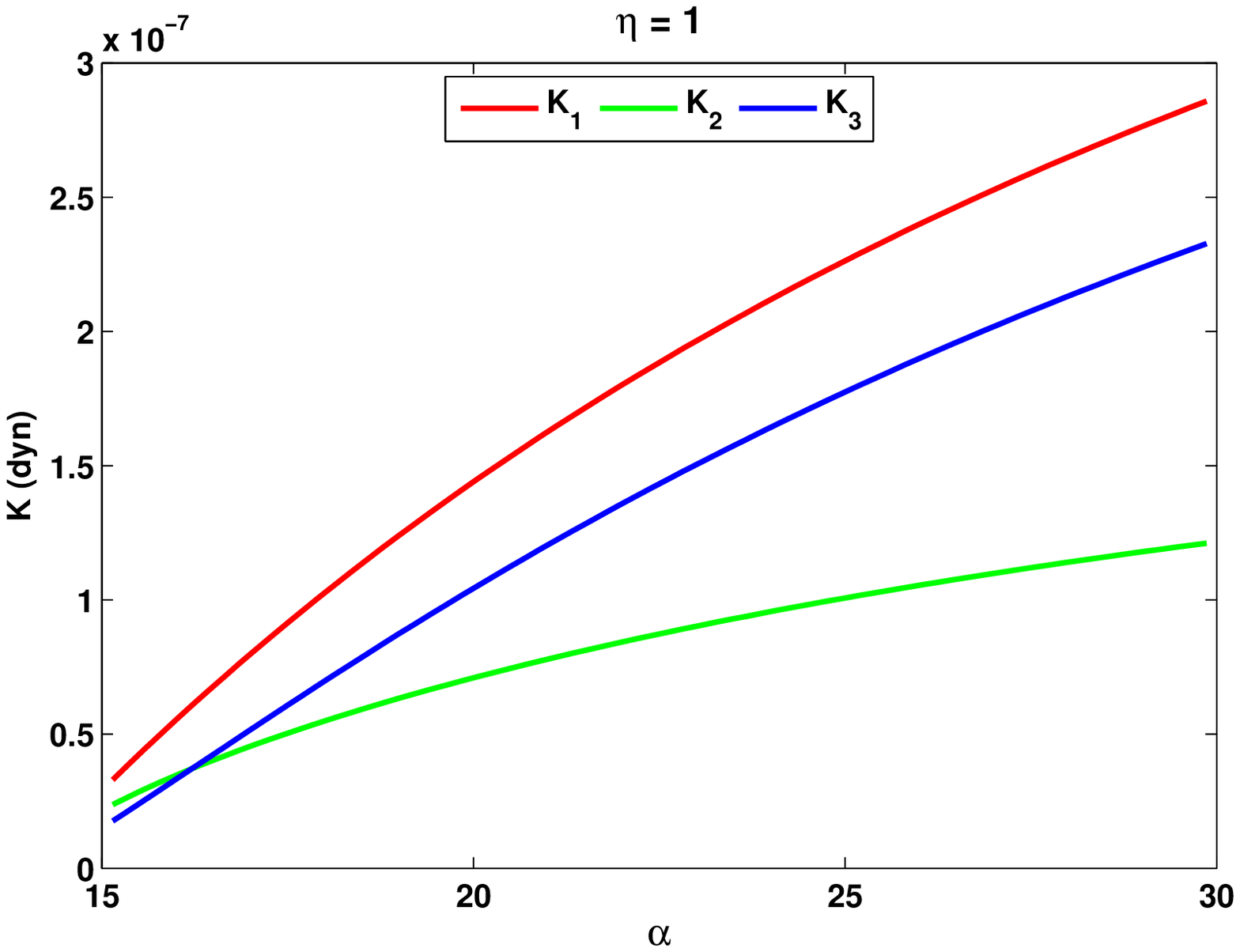}
   \caption{Oseen-Frank elastic constants $K_1, K_2, K_3$ under hard-core potential}
   \label{figoscoef}
\end{figure}

Taking the diameter-to-length ratio $\eta$ as $0.1,~0.3,~0.6$ and $1.0$
respectively, we draw the value curves of the corresponding
Oseen-Frank elastic constants with respect to $\alpha$ in Fig. \ref{figoscoef}.

As the parameters are not accurate and several approximations have
been made in modeling, our elastic constants $K_1, K_2, K_3$ are
close to the physical observation values but not precise.
Nevertheless, we can still conclude some points from them:
\begin{itemize}
\item $ K_1, K_2, K_3$ are keeping positive in value regardless of $\alpha$ and $\eta$, which bounds the distortion energy from below;

\item $K_1, K_2, K_3$ converge in most of the cases as the concentration
factor $\alpha$ gets smaller. It illustrates why we might employ the one-constant
Landau-de Gennes model (in which $K_1=K_2=K_3$) for diluted LC;

\item When the potential density $\alpha$ goes larger, $K_3$ becomes
significantly larger than the other two elastic constants while $K_1$
and $K_2$ are almost remaining at the same level. Such phenomenon have
already been observed in many physical experiments, for instance, the
Frederiks experiment \cite{FZ33};

\item For a typically long rod-like LC molecule, it always has the
following relation:
    $$K_3 \ge K_1 \ge K_2,$$
    no matter how the concentration changes.  {Furthermore, this relation holds
    in the limiting case in which the molecular diameter goes to 0}. This property is also
    supported by experimental results.

\item {When the rod-like LC molecule gets shorter in its shape, elastic
constant $K_1$ will eventually overcome $K_3$ in the low concentration area. In fact, the diameter-to-length ratio of the rod-like
LC molecules are generally observed within the range of $1:4$ to
$1:15$. In this light, this situation might indicate that rod-like molecules which are ``too short" in its length can not be observed in nematic phase.}
\end{itemize}
In addition to the last point, we guess that when the rod-like LC molecule
goes "too long", the nematic phase might also lose its stability to the
smectic phase. We think it will lead us to another fascinating story.

\subsection{The Ericksen's vector model}
While the head-tail symmetry of the LC molecules is respected by both
he molecular model and tensor model, the equivalence of orientation
$\textbf{n}$ and $-\textbf{n}$ was not fully recognized by Oseen.
Consequently, this drawback results in its inability to account for
the complicated line and surface defects that are physically observed.
Particularly in certain circumstances, the Oseen-Frank model will lead
to nonphysical solutions, and even fake ``defects". Besides, as we can
see from the relationship between the order parameter $s$ and
the concentration factor $\alpha$ in the homogeneous case, $s$ moves
rapidly when there is a slight change of $\alpha$ in the nematic phase
area. In other words, while it is
suitable to treat $\alpha$ (or $c$ equivalently) as a constant in nematic
modeling, it might not be proper to set the scalar order parameter $s$ invariant.

Ericksen \cite{Eri91} extended the Oseen-Frank model by relaxing the assumption
of a spatially invariant degree of orientational order $s$. We can also
derive the Ericksen model by imposing uniaxial constraint in our new
tensor model. The only difference is that we have to keep terms containing
$s$ or gradient of $s$.
 Unlike the situation in deriving the Oseen-Frank model, the bulk energy
 denoted in the Ericksen model is not an additive constant, but a
 function of $s$. One can obtain that $F_\text{bulk}^{(Eri)}(s)$
 can be written as follows:
$$F_\text{bulk}^{(Eri)}(s) = k_B T\bigg( c(B_Q:Q-\mathrm{ln}Z)-\frac23 c^2L^2 D s^2\bigg),$$
where both the value of the inner tensor product $B_Q:Q$ and the
value of the normalization constant $Z$ depend only $s$. It is not
difficult to prove that the consistent condition proposed by Ericksen:
$$\lim\limits_{s \to 0} F_\text{bulk}^{(Eri)}(s) = O(s^2) $$
is satisfied here. Besides, we have:
$$F_\text{bulk}^{(Eri)}(s) \to \infty, \qquad \mathrm{if}\quad s \to 1_{-} \quad \mathrm{or}\quad s \to -\frac12_{+} .$$

The Ericksen elastic energy density truncated to the order of $s^2$ reads:
\begin{align}\nonumber
2F_\text{elastic}^{(Eri)} = & K_1 s^2(\mathrm{div} \textbf{n})^2 + K_2 s^2
(\textbf{n}\cdot (\nabla \times \textbf{n}))^2 + K_3 s^2|\textbf{n} \times (\nabla\times \textbf{n})|^2\\
& + (K_2+K_4) s^2(\mathrm{tr}(\nabla \textbf{n})^2-(\mathrm{div} \textbf{n})^2)\nonumber\\
& + l_1|\nabla s|^2 + l_2 |\nabla s \cdot \textbf{n}|^2 + l_3 s \nabla
s \cdot \textbf{n} \mathrm{div} \textbf{n} + l_4 s \nabla s \cdot (\nabla \textbf{nn}),
\end{align}
where the elastic constants $K_1, K_2, K_3$ are the same as what we
have computed for the Oseen-Frank model. The rest of the coefficients
can also be derived from our tensor model as well. We omit their exact
expression here for brevity.

%

\section{Modeling for simple smectic liquid crystals}
Smectic liquid crystals are characterized by both orientational and translational
ordering of anisotropic molecules. In simple smectic phases, namely, smectic-A ($\mathrm{S}_\mathrm{A}$)
and smectic-C ($\mathrm{S}_\mathrm{C}$), the translational ordering is one dimensional, and the director of primary
molecular {axis} is either parallel or tilted with respect to the direction of the
wave vector of the corresponding density wave. Characterizing the nematic to
smectic-A phase transition has long been a principal problem in the physics
of liquid crystals. The main complexity of the Nematic-$\mathrm{S}_\mathrm{A}$ transition arises from an
intrinsic coupling between local number density and order parameter. Despite
considerable literature on modeling of  {smectic phase}, there are still unresolved issues due to this complexity. Based on
the whole idea described above, here  {a simple} model can be constructed
to characterize the nematic and smectic-A phase universally.

\subsection{A tensor model for smectic-A phase}
The key point in modeling the smectic phase is to build the layer structure. Compared with the nematic modeling, an additional positional order parameter
must be introduced to describe the modulation of the concentration.
The derivation procedure is almost the same as our nematic modeling except for that the number density parameter $c$ is no longer a spatially invariant constant.
Therefore to assure the whole free energy bounded from below, we need to truncate the elastic energy to the fourth moment of the kernel function:
$$M^{(4)} : = \int G(\textbf{r},\textbf{m},\textbf{m}')\textbf{r}\textbf{r}\textbf{r}\textbf{r}\,\mathrm{d}\textbf{r}.$$
In the case of hard-core potential, the fourth moment under the original coordinate can be written as (see Appendix for details):
\begin{align}
M^{(4)} &=R_1(\um,\um')(\udeltayier\udeltasansi)_{\mathrm{sym}}
+R_2(\um,\um')(\udeltayier\um\um+\udeltayier\um'\um')_{\mathrm{sym}} \nonumber \\
&\quad\quad+R_3(\um,\um')(\udeltayier\um\um')_{\mathrm{sym}}
+R_4(\um,\um')(\um\um\um\um+\um'\um'\um'\um')\nonumber \\
&\quad\quad+R_5(\um,\um')(\um\um\um'\um')_{\mathrm{sym}}
+R_6(\um,\um')(\um\um\um\um'+\um'\um'\um'\um)_{\mathrm{sym}},
\label{fourthmoment}
\end{align}
where
$$
\left\{
\begin{array}{ll}
{R_1(\um,\um')=L^6D\Big(\frac{2\sin{\gamma}~\eta^4}{15}+\frac{\pi\eta^5}{12}+\frac{4\pi\eta^6}{105}\Big)}, \\
{R_2(\um,\um')=L^6D\Big(\frac{\sin{\gamma}~\eta^2}{18}+\frac{\pi\eta^3}{12}+\frac{\pi\eta^4}{15}
+\frac{\pi\eta^5}{24}+\frac{2\eta^4}{15}\frac{1}{\sin{\gamma}}\Big)}, \\
{R_3(\um,\um')=L^6D\Big(\frac{(\pi-2\gamma)\eta^4}{15}-\frac{2\eta^4}{15}\frac{\cos{\gamma}}{\sin{\gamma}}\Big),}
 \\
{R_4(\um,\um')=L^6D\Big\{\Big(\frac{\sin{\gamma}}{40}+\frac{3\pi\eta}{40}+\frac{\pi\eta^2}{12}+\frac{\pi\eta^3}{
8}-\frac{\pi \eta^5}{24}\Big)+\frac{\eta^3}{3}\frac{1}{\sin{\gamma}}
    -\frac{2\eta^4}{15}\frac{1}{\sin^3{\gamma}}\Big\}},  \\
{R_5(\um,\um')=L^6D\Big\{\Big(\frac{\sin{\gamma}}{72}+\frac{\pi\eta}{24}+\frac{\pi\eta^2}{12}+\frac{\pi\eta^3}{
8}\Big)
    +\Big(\frac{\eta^2}{9}+\frac{2\eta^4}{15}\Big)\frac{1}{\sin{\gamma}}
    -\frac{2\eta^4}{15}\frac{\cos^2{\gamma}}{\sin^3{\gamma}}\Big\}},  \\
{R_6(\um,\um')=L^6D\Big\{\frac{(\pi-2\gamma)\eta^2}{12}
    -\frac{\eta^2}{6}\frac{\cos{\gamma}}{\sin{\gamma}}
    +\frac{2\eta^4}{15}\frac{\cos^3{\gamma}}{\sin^3{\gamma}}\Big\}.}
\end{array}
\right.
$$
For simplicity we just expand the fourth moment to $O(\eta)$:
\begin{align*}
M^{(4)}\approx & \pi L^6D\Big[[\mu_{11}(\mm\mm\mm\mm+\mm'\mm'\mm'\mm')+(\mu_{21}+\mu_{22}P_2(\mm\cdot\mm')^2)(\mm\mm\mm'\mm')
_{\mathrm{sym}}\Big],
\end{align*}
where
\begin{equation*}
\mu_{11}=\frac{1}{160}+\frac{3\eta}{40},\quad\quad\quad \mu_{21}=\frac{1}{288}+\frac{\eta}{24},\quad\quad\quad \mu_{22}=-\frac{5}{2304.}
\end{equation*}
Thus the fourth order elastic energy reads
\begin{align*}
F^{(4)}_{\text{elastic}}=&\int_{\Omega}\int_{\BS}\int_{\BS}\int_{\Omega} f(\xx,\mm)G(\mm,\mm',\rr)(r_i\partial_i)^4
\big\{f(\xx,\mm')\big\}~\ud\rr\ud\mm'\ud\mm\ud\xx\nonumber\\
\approx &\frac{\pi L^7\eta k_BT}{24}\int_{\Omega}\int_{\BS}\int_{\BS} \bigg\{ 2\mu_{11}~f(\xx,\mm)~m_im_jm_km_l~\partial_{ijkl}\big\{f(\xx,\mm')\big\}\ud\mm'\ud\mm  \nonumber\\
&\quad\quad\quad+(\mu_{21}-\frac{1}{2}\mu_{22})f(\xx,\mm)(m_im_jm'_km'_l)_{\mathrm{sym}}~\partial_{ijkl}\big\{f(\xx,\mm')\big\}
\ud\mm'\ud\mm \nonumber\\
&\quad\quad\quad+\frac{3}{2}\mu_{22}f(\xx,\mm)(\mm\cdot\mm')^2(m_im_jm'_km'_l)_{\mathrm{sym}}~\partial_{ijkl}\big\{f(\xx,\mm')\big\}
\ud\mm'\ud\mm\bigg\} \ud\xx.
\end{align*}

$F^{(4)}_{\text{elastic}}$ can be similarily written in the Q-tensor form (see Appendix for details). Since we introduce $Q_4$ in $F^{(4)}_{\text{elastic}}$, we
also truncate the approximation of $|\textbf{m}\times\textbf{m}'|$ at the fourth order Legendre polynomial.
Now the total free energy functional is readily given as
$$F_{\text{total}}=F_{\text{bulk}}+F^{(2)}_{\text{elastic}}+F^{(4)}_{\text{elastic}}.$$
After proper substitution $\xx=\xx/L$, we finally reach a dimensionless Q-tensor model with the free energy functional as follows:
\begin{align}
&~~~~F[c(\xx),Q_2(\xx)] \nonumber\\
&=\int_{\Omega} c(\ln{c}+B_Q:Q_2-\ln{Z})\ud\xx+\frac{\alpha}{2}\int_{\Omega}\big\{ E_{11}c^2+E_{12}|cQ_2|^2+E_{13}|cQ_4|^2 \nonumber\\
&\quad+E_{21}|\nabla{c}|^2+E_{22}|\nabla(cQ_2)|^2+E_{23}|\nabla(cQ_4)|^2+E_{24}~\partial_i(cQ_{2ij})\partial_j(c) \nonumber \\*[0.2cm]
&\quad
+E_{25}~\partial_i(cQ_{ik})\partial_j(cQ_{jk})+E_{26}~\partial_i(cQ_{4ijkl})\partial_j(cQ_{2kl})
+E_{27}~\partial_i(cQ_{4iklm})\partial_j(cQ_{4jklm}) \nonumber \\*[0.2cm]
&\quad+E_{31}|\nabla^2{c}|^2+E_{32}~\partial_{ij}(cQ_{2pq})\partial_{ij}(cQ_{2pq})+E_{33}~\partial_{ij}(cQ_{2ij})\partial_{kl}(cQ_{2kl})
\nonumber\\*[0.2cm]
&\quad+E_{34}~\partial_{ik}(cQ_{2ip})\partial_{jk}(cQ_{2jp})+E_{35}~\partial_{ij}(cQ_{2ij})\partial_{kk}(c)+E_{36}~\partial_{ij}(cQ_{
4ijkl})\partial_{kl}(c) \nonumber\\*[0.2cm]
&\quad+E_{37}~\partial_{ij}(cQ_{4ijpq})\partial_{kk}(cQ_{2pq})+E_{38}~\partial_{ij}(cQ_{4ijkp})\partial_{kl}(cQ_{2lp}) \nonumber \\*[0.2cm]
&\quad+E_{39}~\partial_{ij}(cQ_{4ijpq})\partial_{kl}(cQ_{4klpq})\big\}\ud\xx.
\label{Q4smectic}
\end{align}
Here $E_{ij}$ only depend on $\eta$. $B_Q, Z, Q_4$ are determined by $Q_2$ with Bingham closure.
This model has only two dimensionless parameters: $\eta=D/L$ and $\alpha=\pi c_0L^2D$ where
$$\int c(\xx)\ud\xx = c_0.$$ Notice that $c(\xx)$ now satisfies the constraint
$$\frac{1}{|\Omega|}\int_{\Omega} c(\xx) \ud\xx =1.$$

\subsection{One-dimensional model and numerical results}
Based on the tensor model derived above, we consider a one-dimensional model with following two assumptions:
{\begin{itemize}
\item{} $f(\xx,\mm)$
only depends on $x$-axial component and is a periodic
function with period $d$;
\item{} Orientation distribution are uniaxial and the director $\nn$ is a constant parallelling with $x$ axis.
\end{itemize}}
Notice all non-trivial tensor components needed are given as
\begin{align*}
&Q_{11}=\frac{2}{3}S_2, ~~~~\quad\quad\quad Q_{22}=Q_{33}=\frac{1}{3}S_2, \\
&Q_{1111}=\frac{8}{35}S_4, \quad\quad\quad Q_{1122}=Q_{1133}=-\frac{4}{35}S_4.
\end{align*}
Here $S_4(x)$ is determined by Bingham Distribution or say $r(x)$:
\begin{align*}
S_4(x)&=\frac{35\displaystyle\int_0^1t^4\exp{(r(x)t^2)}~\ud t}{8\displaystyle\int_0^1\exp{(r(x)t^2)}~\ud
t}-\frac{5}{2}S_2(x)-\frac{7}{8}.
\end{align*}
Now (\ref{Q4smectic}) can be furthermore reduced to a one-dimension model with free energy functional as
follows:
\begin{align}
&\bar{F}_{d}[c(x),S_2(x)] \nonumber \\
\triangleq &\int_0^d c(x)\Big(\ln c(x)+\frac{2}{3}r(x)S_2(x)-\ln Z(x)\Big)\ud x\nonumber \\
&+\frac{\alpha}{2}\int_0^d\bigg\{ \big(N_{11}+N_{12}S_2^2(x)+N_{13}S_4^2(x)\big)c^2(x) \nonumber \\
&\quad\quad~~~-N_{21}
\Big(\udxyi\big(c(x)\big)\Big)^2-N_{22}\Big(\udxyi\big(c(x)S_2(x)\big)\Big)^2-N_{23}(\udxyi\big(c(x)S_4(x)\big)\Big)^2\nonumber\\
&\quad\quad~~~-N_{24}\udxyi\big(c(x)S_2(x)\big)\udxyi\big(c(x)\big)-N_{25}\udxyi\big(c(x)S_2(x)\big)\udxyi\big(c(x)S_4(x)\big)
\nonumber\\
&\quad\quad~~~+N_{31}
\Big(\udxer\big(c(x)\big)\Big)^2+N_{32}\Big(\udxer\big(c(x)S_2(x)\big)\Big)^2\bigg\}\ud x.
\label{1dsmectic}
\end{align}
All the coefficients $N_{ij}$ could be calculated from the moment calculation and expansion above:
\begin{align*}
&N_{11}=\frac{1}{2}+2\eta+\frac{4\eta^2}{3}, \quad\quad\quad\quad\quad N_{12}=-\frac{5}{16},\quad\quad\quad\quad\quad\quad~
N_{13}=-\frac{9}{128}, \\
&N_{21}=\frac{1}{72}+\frac{\eta}{9}+\frac{5\eta^2}{18}+\frac{\eta^3}{3}+\frac{2\eta^4}{15},
\quad\quad\quad N_{22}=-\frac{55}{4032}-\frac{(432\ln{2}-367)\eta^2}{4704},\\
&N_{23}=\frac{(365-2048\ln{2})\eta^2}{12544},\quad\quad\quad\quad\quad\quad
N_{24}=\frac{7}{288}+\frac{2\eta}{9}+\frac{7\eta^2}{18}+\frac{\eta^3}{6},\\
&N_{45}=-\frac{1}{112}-\frac{(12\ln{2}-10)\eta^2}{49},\\
&N_{31}=\frac{11}{57600}+\frac{13\eta}{5400},\quad\quad\quad\quad\quad\quad\quad\quad~N_{32}=\frac{107}{451584}+\frac{\eta}{216}.
\end{align*}
In (\ref{1dsmectic}) we just keep two second derivative terms and drop
other four terms to make our model as simple as possible. In the derivation of
the above coefficients, we have made several approximations and truncations.
As a result, the numerical value might not be accurate and
perhaps some information such as attraction effect and temperature
dependence is lost. Then it might cause the free energy functional
not bounded from below. {Actually, the optimal solution depends on these
two coefficients sensitively. Therefore, to find physical solutions,
we modify $N_{31}$ and $N_{32}$ in a reasonable range without changing their orders.}
 The main point is that (\ref{1dsmectic}) should be
an effective energy form to capture the smectic-A phase.

To solve the {optimization problem}
\begin{eqnarray*}
\min \limits_{c(x), S_2(x), d}~\bigg\{\frac{\bar F_{d}}{d}\bigg\}, \qquad
s.t.\quad \frac{1}{d}\int_0^d c(x) \ud x =1,
\end{eqnarray*}
we use the spectral method. For this, we have to expand $c(x), c(x)S_2(x)$ in terms of fourier bases:
\begin{equation*}
1, \cos{\frac{2\pi}{d}x}, \sin{\frac{2\pi}{d}x}, \cdots,
\end{equation*}
and truncate $c(x), c(x)S_2(x), c(x)S_4(x)$ at order $n_1, n_2, n_3$, $i.e.$
\begin{align*}
c(x)&=1+\sum_{n=1}^{n_1}u_n\cos{\frac{2n\pi}{d}x},\\
c(x)S_2(x)&=\sum_{n=0}^{n_2}v_n\cos{\frac{2n\pi}{d}x},\\
c(x)S_4(x)&=\sum_{n=0}^{n_3}w_n\cos{\frac{2n\pi}{d}x},\\
\ln c(x)+\frac{2}{3}r(x)S_2(x)&-\ln Z(x)=\sum_{n=0}^{n_1}t_n\cos{\frac{2n\pi}{d}x}.
\end{align*}
Notice that $u_n, v_n$ are variables and $w_n, t_n$ are determined {by solving $r(x)$
and integration with FFT. The minimum of the free energy can be found
by standard method, for example, the steepest descent method.}

\begin{figure}[!h]
\centering \mbox{
\subfigure[One dimensional phase diagram for isotropic, nematic, smectic-A phase.]{\includegraphics[width=7.4cm,height=5cm]{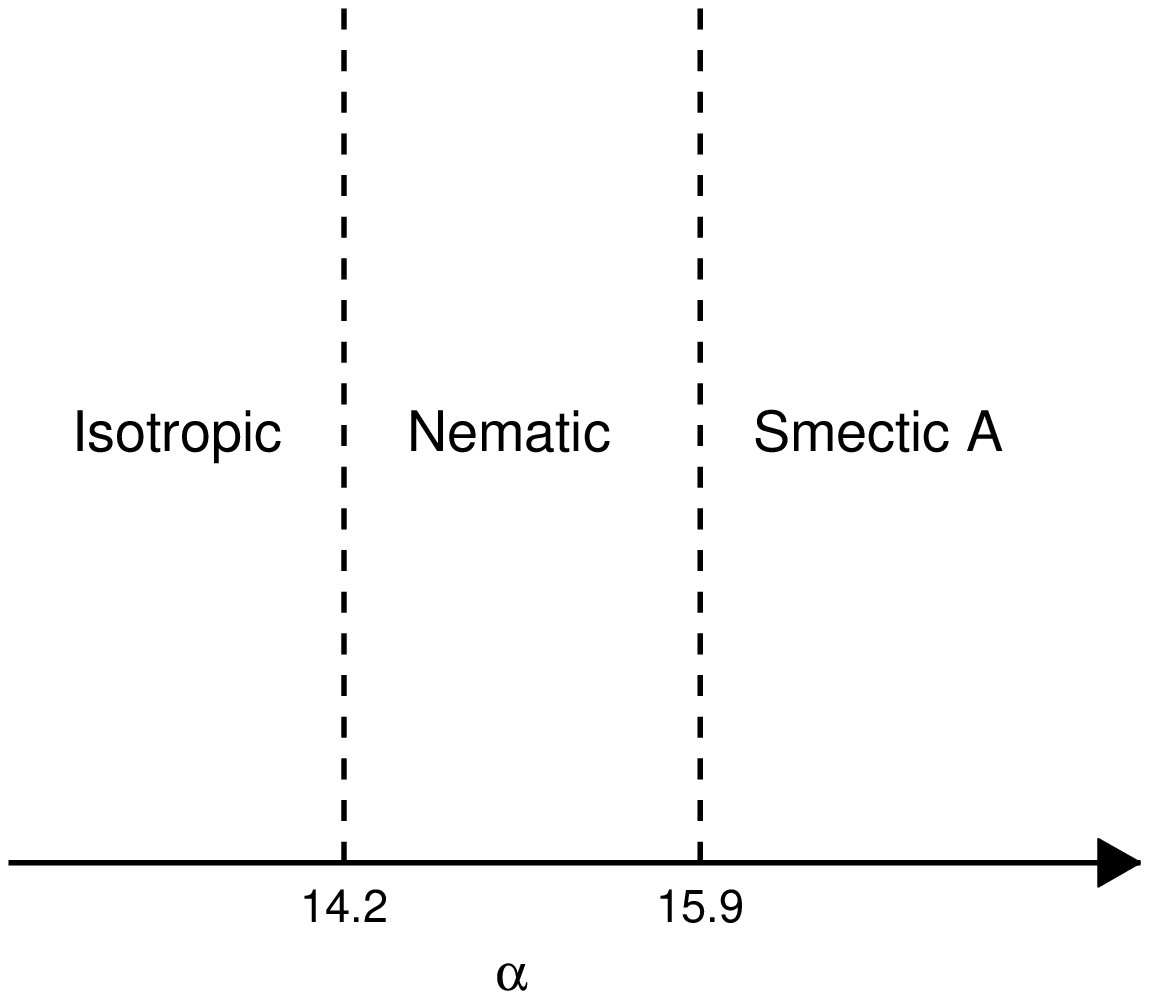}}
\quad\quad
\subfigure[Profiles of the local number density $c(x)$
(dashed) and local nematic order parameter (solid) in the smectic-A
phase at $\alpha=17$: $d = 1.519~L$.]
{\includegraphics[width=7.2cm,height=5.1cm]{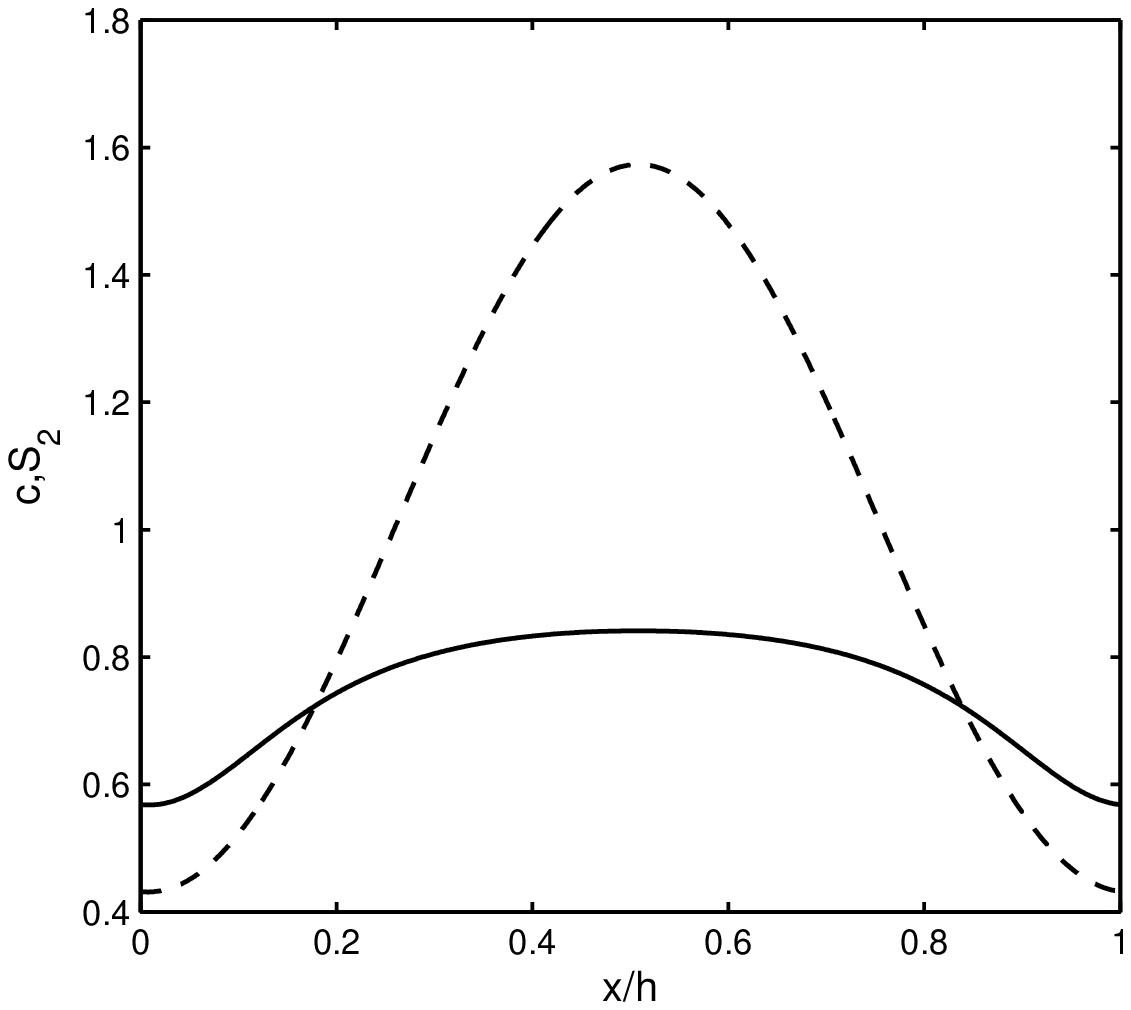}}
} \caption{Typical phase diagrams, local number density and order parameter profiles. $\eta=1/10,~N_{31}=N_{32}=0.00089$.}
\label{fig:phase1}
\end{figure}

Fig. \ref{fig:phase1} (a) presents a typical phase diagram of  {three phases}. Nematic phase loses stability  {as concentration increases}. The smectic layer periodicity $d$ increases slightly with the increasing of concentration and spans from $1.516~L$ to $1.532~L$, which is quite reasonable according to existing experimental results \cite{dGP95}. It is also worth noting that the boundary of nematic phase in the phase diagram depends on $N_{31}, N_{32}$ and in some cases the system only undergoes direct isotropic-$\mathrm{S}_\mathrm{A}$ phase transition, which agrees with some experimental results. Fig. \ref{fig:phase1} (b) reveals typical fluctuations of local number density and nematic order parameter.

\par

\subsection{Discussion on smectics modeling}
Nematic-$\mathrm{S}_\mathrm{A}$ phase transition has been addressed
theoretically via different approximations for a long time. Here we briefly review two popular smectic models.

McMillan \cite{McM71} first put forward a specific model, which is an extension
of the Maier-Saupe mean-field theory, characterizing the smectic-A phase by a
density modulation. Applying Landau expansions for the entropy term, this model
describes a continuous nematic-$\mathrm{S}_\mathrm{A}$ transition and predicts the existence
of a triple critical point where the nematic, isotropic, and smectic-A phase
meet.
Besides the orientational order parameter defined in the Maier-Saupe theory, an order parameter
that describe the positional order of the LC is introduced in McMillan's model.
Since the smectic-A phase is uniaxial, the positional order parameter is given by:
$$\sigma_p = \langle \cos (\frac{2\pi z}{d}) (\frac32 \cos^2 \theta -\frac12) \rangle $$
where $z$ and $d$ are the position of the molecule and the layer thickness. The postulated interaction potential reads:
$$U(\theta, z)=-U_1(S+\alpha \sigma_p (\frac{2\pi z}{d}))(\frac32 \cos^2 \theta - \frac12) $$
where the constant $\alpha$ refers to the strength between adjacent molecules.

Another successful model for smectics is the
{Chen-Lubensky model \cite{CL76}}, which is defined in terms of the director field $\textbf{n}$ and
the complex valued smectic order parameter
$$\Psi(\textbf{x})=\rho(\textbf{x})e^{i\phi(\textbf{x})}, $$
where $\phi$ refers to the layers and consequently $\nabla \phi$ is the
direction of the layer normal. For a perfect nematic phase $\Psi=0$ while for a smectic phase $\Psi$ will take on complex values.

It is a phenomenological vector model based on the Landau-Ginzburg mean-field
theory. The free energy consists of two parts:
$$F^{(CL)}=\int_{\Omega}F_S + F_N \,\mathrm{d}\textbf{x},  $$
where $F_N=F^{(OF)}$ is the Oseen-Frank energy density for a nematic. The smectic free energy density $F_S$ could be designed to describe either nematic-$\mathrm{S}_\mathrm{A}$ (or A*) phase transition or nematic-smectic-C (or C*) phase transition. If we only hope to model smectic-A (or A*) phase, the free
energy can be simplified to the following form of Landau-de Gennes energy density:
\begin{align}
F_S=a(T-T_{N\mathrm{S}_\mathrm{A}})|\Psi|^2 +\frac{1}{2}g|\Psi|^4 - C_{\perp}|\textbf{D}\Psi|^2.
\label{energy:sm_chen}
\end{align}
In (\ref{energy:sm_chen}) we have: $\textbf{D}=\nabla - iq \textbf{n}$; $q \sim 1/d$ where $d$
is the smectic layer thickness; $T_{N\mathrm{S}_\mathrm{A}}$ is the temperature where
nematic-$\mathrm{S}_\mathrm{A}$ phase transition occurs; the coefficients $q$,$g$ are positive and $C_{\perp}\le 0$ is necessarily required.

The McMillan model and Chen-Lubensky model both require a prior knowledge of the
layer thickness $d$ before modeling. Furthermore, the Chen-Lubensky model lacks clear physical interpretations of its various
coefficients especially for $C_{\perp}$ which determines the type
of smectic phase. The McMillan model does not take the smectic-C phase into consideration.
{In Chen-Lubensky model, axially-symmetry like a `tilt
uniaxial'  is assumed to describe smectic-C phase.} However, physical observations
show that LC is biaxial in the smectic-C phase.

Compared with the above two models, layer thickness $d$ need not a priori
in our Q-tensor model for smectic-A phase since it is obtained by minimizing the free energy.
 Although the dependence of coefficient of higher orders on temperature should be studied more deeply,
 all coefficients in our model could be determined by physical parameters,
 which provides opportunities to further testify our model through experimental
 results. We also point out that the Q-tensor model might also be able to describe
 the smectic-C phase if we relax the assumption to allow director $\nn$ tilted
 with respect to the direction of wave vector and introduce biaxial
 approximation. All these works would lead to a subject of a separate study.

\section{Summary}
 {We have proposed a multiscale modeling hierarchy for liquid crystals connecting Onsager's molecular theory, $Q$-tensor theory and Oseen-Frank theory explicitly.} As an important example, we discussed the simple
stick-shape molecule with head-tail symmetry.
{Once the intermolecular potential is decided (such as hard-core
 potential or Lennard-Jones potential)}, we are able to
write the integral form of the molecular model.   {Employing Taylor
expansion,} the differential form of the molecular model, and
$Q$-tensor can be introduced to simplify the expression. Different
truncations and approximations for the high order moments of the
kernel function will lead to different models.

To model nematic phase, it is sufficient to truncate at the second
order of derivatives. Then we obtain a new $Q$-tensor type model
for nematic liquid crystals, which not only remains sensitive to the
macroscopic properties, but also  {take into account} the molecular structure
and mechanics as well. Its distortion energy is bounded from below,
and the existence of a physically meaningful minimizer is guaranteed.
The coefficients in this model are
entirely determined by molecular parameters. In addition, the order
parameter satisfies the physical constraint on eigenvalues naturally.
From this model, we can also recover the  {Oseen-Frank energy coefficients} by the molecular parameters.
%

On the other hand, if we truncate at the fourth order of derivatives,
the obtained model can describe the smectic phase, which is convinced by
some  {numerical results.}

Several new models for a  {variety of} molecular
structures and interactions can be derived by following the same
procedure. All these models shall meet the physical constraints
and all the constants shall be interpretable and easy to
determine. We believe that these models will help to study other different LC phases.

\section{Appendix}
\subsection{High order traceless symmetric tensor}
For any  axisymmetric function $f(\mm)=f(\mm\cdot\nn)$, we want to calculate
$$\langle\underbrace{\mm\otimes\mm\otimes\cdots\otimes\mm}_{k\text{ times}} \rangle_{f}.$$
This motivates us to introduce the $k$-order traceless symmetric tensor on the unit sphere.
We use $\widehat{1},\widehat{2},\cdots\in\{1,2,3\}$ to denote subscripts. For $\mm\in\BS$, we define the $(k+2l)$-order symmetric tensor as follows
\begin{align}
\sigma(\mm;k,2l)_{\widehat{1}\widehat{2}\cdots\widehat{m+2l}}=
\big(m_{\widehat{1}}m_{\widehat{2}}\cdots{m}_{\widehat{k}}
\delta_{\widehat{k+1}\widehat{k+2}}\cdots\delta_{\widehat{k+2l-1}\widehat{k+2l}}\big)_{\text{sym}},
\end{align}
where $(\cdot)_{\text{sym}}$ means the symmetrization of the tensor. For example,
\begin{align*}
\sigma(\mm, 1,2)&=m_\alpha\delta_{\beta\gamma}+m_\beta\delta_{\alpha\gamma}+m_\gamma\delta_{\alpha\beta},\\
\sigma(\mm, 2,2)&=m_{\alpha}m_{\beta}\delta_{\gamma\mu}+
m_{\gamma}m_{\mu}\delta_{\alpha\beta}+m_{\alpha}m_{\gamma}\delta_{\beta\mu}
+m_{\beta}m_{\mu}\delta_{\alpha\gamma}+m_{\alpha}m_{\mu}\delta_{\beta\gamma}+
m_{\beta}m_{\gamma}\delta_{\alpha\mu},\\
\sigma(\mm,0,4)&=\delta_{\alpha\beta}\delta_{\gamma\mu}
+\delta_{\alpha\gamma}\delta_{\beta\mu}+\delta_{\alpha\mu}\delta_{\beta\gamma}.
\end{align*}
Direct computation shows that $\sigma(k,2l)$ is a sum of $\frac{(k+2l)!}{k!l!2^l}$ different tensors.
If we contract any two subscripts of $\sigma(k,2l)$, we  obtain a
sum of some $(k+2l-2)$-order symmetric tensors, i.e,
\begin{align*}
\text{Contract}\big[\sigma(k,2l)\big]=&3\sigma(k,2l-2)+\sigma(k-2,2l)+2k\sigma(k,2l-2)+2(l-1)\sigma(k,2l-2)\\
=&(2k+2l+1)\sigma(k,2l-2)+\sigma(k-2,2l).
\end{align*}

Let
\begin{align*}
\Xi_n(\mm)=\sum_{l=0}^{[\frac{n}{2}]}a_l\sigma(\mm;n-2l,2l), \qquad a_0=1.
\end{align*}
We need $\Xi$ to be {\it trace free}, that is, if we contract any two
subscripts of $\Xi$, the result tensor should be 0. Therefore,
\begin{align*}
\sum_{l=0}^{[\frac{n}{2}]}a_{n,l}\Big\{(2n-2l+1)\sigma(n-2l,2l-2)+\sigma(n-2l-2,2l)\Big\}=0.
\end{align*}
Hence
\begin{align*}
a_{n,0}\sigma(n-2,0)
+a_{n,1}\big[(2n-1)\sigma(n-2,0)+\sigma(n-4,2)\big]\\
+a_{n,2}\big[(2n-3)\sigma(n-4,2)+\sigma(n-6,4)\big]+\cdots=0.
\end{align*}
Thus we have
\begin{align*}
a_{n,l-1}=-(2n-2l+1)a_{n,l}, \text{ for }l\ge1.
\end{align*}
Therefore
\begin{align*}
\Xi_n(\mm)=&\sigma(n,0)-\frac{1}{2n-1}\sigma(n-2,2)+\frac{1}{(2n-1)(2n-3)}\sigma(n-4,4)\nonumber\\
&-\frac{1}{(2n-1)(2n-3)(2n-5)}\sigma(n-6,6)+\cdots.
\end{align*}
The following are some examples:
\begin{align*}
&\Xi_1(\mm)=\mm;\\
&\Xi_2(\mm)=m_\alpha m_\beta-\frac{1}{3}\delta_{\alpha\beta};\\
&\Xi_3(\mm)=m_{\alpha}m_{\beta}m_{\gamma}-\frac{1}{5}\Big(m_\alpha\delta_{\beta\gamma}
+m_\beta\delta_{\alpha\gamma}+m_\gamma\delta_{\alpha\beta}\Big);\\
&\Xi_4(\mm)={m}_{\alpha}{m}_{\beta}{m}_{\gamma}{m}_{\mu}
-\frac{1}{7}\Big(m_{\alpha}m_{\beta}\delta_{\gamma\mu}+
m_{\gamma}m_{\mu}\delta_{\alpha\beta}+m_{\alpha}m_{\gamma}\delta_{\beta\mu}
+m_{\beta}m_{\mu}\delta_{\alpha\gamma}\nonumber\\
&\qquad\qquad+m_{\alpha}m_{\mu}\delta_{\beta\gamma}+
m_{\beta}m_{\gamma}\delta_{\alpha\mu}\Big)+\frac{1}{35}\Big(\delta_{\alpha\beta}\delta_{\gamma\mu}
+\delta_{\alpha\gamma}\delta_{\beta\mu}+\delta_{\alpha\mu}\delta_{\beta\gamma}\Big);\\
&\Xi_5(\mm)={m}_{\alpha}{m}_{\beta}{m}_{\gamma}{m}_{\mu}m_\nu-\frac{1}{9}\Big(m_{\alpha}m_{\beta}m_{\gamma}
\delta_{\mu\nu}+\cdots\Big)+\frac{1}{63}\Big(m_\alpha\delta_{\beta\gamma}\delta_{\mu\nu}+\cdots\Big).
\end{align*}

Let $P_n$ be the Legendre polynomials:
\begin{align}
P_n(x)=2^n\sum_{k=0}^nx^k{n \choose k}{\frac{n+k-1}{2} \choose n}=
\frac{1}{n!}\sum_{k=0}^{[\frac{n}{2}]}(-1)^k{(2n-2k-1)!!}
\frac{n!}{(n-2k)!k!2^k}x^{n-2k}.
\end{align}
The following proposition enables us to calculate
$\langle\underbrace{\mm\otimes\mm\otimes\cdots\otimes\mm}_{k\text{ times}} \rangle_{f}$
for any  axisymmetric function $f(\mm)=f(\mm\cdot\nn)$.
\begin{proposition}
For any axisymmetric function on $\BS$: $f(\mm)=f(\mm\cdot\nn)$, define
$S_k[f]=\int_\BS P_k(\mm\cdot\nn)f(\mm)\ud\mm$. Then we have
\begin{align}
\int_\BS\Xi_k(\mm,n)f(\mm)\ud\mm=S_k(f)\Xi_k(\nn).
\end{align}
\end{proposition}

\subsection{The calculation of the zero-th, second and fourth moment for the hard-core interaction potential}
Recall that we separate the whole area into the following three regions:
\begin{itemize}
\item region A (body-body): a $2D$-high parallelepiped whose section is a rhombus with side-length $L$ and angle $\gamma$;
\item region B (body-end): four semi-columns with side-length $L$ and radius $D$;
\item region C (end-end): four radius $D$ sphere at the corner.
\end{itemize}

For the zero-th moment, in region A let
$$(x,y,z) \to ((u+v)\cos \beta, (u-v)\sin \beta, z),$$
$$u\in [-L/2,L/2],\qquad v \in [-L/2,L/2],\qquad z\in [-D,D],$$
then we have
$$\left|\frac{\partial(x,y,z)}{\partial(u, v,z)}\right| = \sin \gamma.$$

In region B, let
$$(x,y,z) \to (t\cos \beta + r \sin\beta \sin\theta, (L-t)\sin\beta + r\cos\beta\sin\theta, r\cos \theta),$$
$$t\in[0,L],\qquad r\in[0,D], \qquad \theta\in [0,\pi],$$
then we have
$$\left|\frac{\partial(x,y,z)}{\partial(t,r,\theta)}\right| = r.$$

In region C, on the one hand, define the region that crosses the real part of $y$-axis as C\Rmnum{1}. Using coordinates which satisfies:
$$(x,y,z) = (r\sin\theta \cos \phi, L\sin\beta + r\sin\theta \sin \phi, r\cos\theta),$$
$$r\in[0,D],\qquad \theta\in [0,\pi], \qquad \phi \in[\frac{\pi}{2}-\beta,\frac{\pi}{2}+\beta],$$
then we have
$$\left|\frac{\partial(x,y,z)}{\partial(r,\theta,\phi)}\right| =r^2 \sin \theta.$$

On the other hand, define the region that crosses the real part of $x$-axis as C\Rmnum{2}, and let
$$(x,y,z) \to (L\cos \beta + r\sin\theta \cos \phi, r\sin\theta \sin\phi, r\cos\theta),$$
$$r\in[0,D],\qquad \theta\in [0,\pi], \qquad \phi \in[-\frac{\pi}{2}+\beta,\frac{\pi}{2}-\beta].$$
We also have:
$$\left|\frac{\partial(x,y,z)}{\partial(r,\theta,\phi)}\right| =r^2 \sin \theta$$

Summing up all the above regions, the total excluded-volume (i.e. the zero moment of the hard-core potential) reads:
$$\int G(|\textbf{r}|,\textbf{m},\textbf{m}') \,\mathrm{d}\textbf{r} = 2L^2D\sin \gamma + 2\pi D^2 L + \frac{4}{3}\pi D^3.$$

Following the same process, we can work out the second moment. In region A, we have
\begin{align*}
\displaystyle{\int_{\text{A}} G(|\ur|,\um,\um')
\left[\begin{array}{c}
    r_1^2 \\ [0.06cm]
    r_2^2\\ [0.06cm]
    r_3^2
\end{array}\right]\ud\ur}
=~&\sin{\gamma}\int_{-D}^{D}\int_{-L/2}^{L/2}\int_{-L/2}^{L/2}
\left[\begin{array}{c}
    (u+v)^2\cos^2{\beta}  \\ [0.06cm]
    (u-v)^2\sin^2{\beta}\\ [0.06cm]
    z^2
\end{array}\right]~\ud u\ud v\ud z \\
=~&\left[\begin{array}{c}
    L^4D\cos^2{\beta}\sin{\gamma}/3  \\ [0.06cm]
    L^4D\sin^2{\beta}\sin{\gamma}/3\\ [0.06cm]
    2L^2D^3\sin{\gamma}/3
\end{array}\right].
\end{align*}
In region B, we can get
\begin{align*}
&\displaystyle{\int_\text{B} G(|\ur|,\um,\um')
\left[\begin{array}{c}
    r_1^2 \\ [0.06cm]
    r_2^2\\ [0.06cm]
    r_3^2
\end{array}\right]\ud\ur} \\
=~&4\int_{0}^{L}\int_{0}^{D}\int_{0}^{\pi}
r\left[\begin{array}{c}
    (t\cos{\beta}+r\sin{\beta}\sin{\theta})^2  \\ [0.06cm]
    ((L-t)\sin{\beta}+r\cos{\beta}\sin{\theta})^2 \\ [0.06cm]
    (r\cos{\theta})^2
\end{array}\right]~\ud \theta\ud r\ud t \\
=~&4\int_{0}^{L}\int_{0}^{D}\int_{0}^{\pi}
\left[\begin{array}{c}
    rt^2\cos^2{\beta}+r^2t\sin{\gamma}\sin{\theta}+r^3\sin^2{\beta}\sin^2{\theta}  \\ [0.06cm]
    r(L-t)^2\cos^2{\beta}+r^2(L-t)\sin{\gamma}\sin{\theta}+r^3\sin^2{\beta}\sin^2{\theta} \\ [0.06cm]
    r^3\cos^2{\theta}
\end{array}\right]~\ud \theta\ud r\ud t \\
=~&\left[\begin{array}{c}
    2\pi L^3D^2\cos^2{\beta}/3+4L^2D^3\sin{\gamma}/3+\pi LD^4\sin^2{\beta}/2  \\ [0.06cm]
    2\pi L^3D^2\sin^2{\beta}/3+4L^2D^3\sin{\gamma}/3+\pi LD^4\cos^2{\beta}/2\\ [0.06cm]
    \pi LD^4/2
\end{array}\right].
\end{align*}

And in region C\Rmnum{1} and C\Rmnum{2}, it reads:
\begin{align*}
&\displaystyle{\int_\text{CI} G(|\ur|,\um,\um')
\left[\begin{array}{c}
    r_1^2 \\ [0.06cm]
    r_2^2\\ [0.06cm]
    r_3^2
\end{array}\right]\ud\ur} \\
=~&\int_{0}^{D}\int_{0}^{\pi}\int_{\pi/2-\beta}^{\pi/2+\beta}
r^2\sin{\theta}\left[\begin{array}{c}
    (r\sin{\theta}\cos{\varphi})^2  \\ [0.06cm]
    (L\sin{\beta}+r\sin{\theta}\sin{\varphi})^2 \\ [0.06cm]
    (r\cos{\theta})^2
\end{array}\right]~\ud \varphi\ud \theta\ud r \\
=~&\int_{0}^{D}\int_{0}^{\pi}\int_{\pi/2-\beta}^{\pi/2+\beta}
\left[\begin{array}{c}
    r^4\sin^3{\theta}\cos^2{\varphi}  \\ [0.06cm]
    L^2r^2\sin^2{\beta}\sin{\theta}+2Lr^3\sin{\beta}\sin^2{\theta}\sin{\varphi}+r^4\sin^3{\theta}\sin^2{\varphi} \\ [0.06cm]
    r^4\cos^2{\theta}\sin{\theta}
\end{array}\right]~\ud \varphi\ud \theta\ud r \\
=~&\left[\begin{array}{c}
    2D^5(\gamma-\sin{\gamma})/15  \\ [0.06cm]
    2L^2D^3\gamma\sin^2{\beta}/3+\pi LD^4\sin^2{\beta}/2+2D^5(\gamma+\sin{\gamma})/15\\ [0.06cm]
    2D^5\gamma/15
\end{array}\right],
\end{align*}

\begin{align*}
&\displaystyle{\int_\text{CII} G(|\ur|,\um,\um')
\left[\begin{array}{c}
    r_1^2 \\ [0.06cm]
    r_2^2\\ [0.06cm]
    r_3^2
\end{array}\right]\ud\ur} \\
=~&\int_{0}^{D}\int_{0}^{\pi}\int_{-\pi/2+\beta}^{\pi/2-\beta}
r^2\sin{\theta}\left[\begin{array}{c}
    (L\cos{\beta}+r\sin{\theta}\cos{\varphi})^2  \\ [0.06cm]
    (r\sin{\theta}\sin{\varphi})^2\\ [0.06cm]
    (r\cos{\theta})^2
\end{array}\right]~\ud \varphi\ud \theta\ud r \\
=~&\int_{0}^{D}\int_{0}^{\pi}\int_{-\pi/2+\beta}^{\pi/2-\beta}
\left[\begin{array}{c}
    L^2r^2\cos^2{\beta}\sin{\theta}+2Lr^3\cos{\beta}\sin^2{\theta}cos{\varphi}+r^4\sin^3{\theta}\cos^2{\varphi}  \\ [0.06cm]
    r^4\sin^3{\theta}\sin^2{\varphi} \\ [0.06cm]
    r^4\cos^2{\theta}\sin{\theta}
\end{array}\right]~\ud \varphi\ud \theta\ud r \\
=~&\left[\begin{array}{c}
    2L^2D^3(\pi-\gamma)\cos^2{\beta}/3+\pi LD^4\cos^2{\gamma}/2+2D^5(\pi-\gamma+\sin{\gamma})/15 \\ [0.06cm]
    2D^5(\pi-\gamma-\sin{\gamma})/15  \\ [0.06cm]
    2D^5(\pi-\gamma)/15
\end{array}\right].
\end{align*}

Summing them up, we can finally get the entire second moment matrix $\mathrm{diag}(M_1,M_2,M_3)$ as follows:
\[
\left[\begin{array}{c}
    \displaystyle{L^4D\left\{\left(\frac{\sin{\gamma}}{3}+\frac{2\pi\eta}{3}+\frac{4(\pi-\gamma)\eta^2}{3}+\pi\eta^3\right)\cos^2{\beta}
    +\left(\frac{4\sin{\gamma}~\eta^2}{3}+\frac{\pi\sin^2{\beta}~\eta^3}{2}+\frac{4\pi\eta^4}{15}
    \right)\right\}} \\ [0.3cm]
    \displaystyle{L^4D\left\{\left(\frac{\sin{\gamma}}{3}+\frac{2\pi\eta}{3}+\frac{4\gamma\eta^2}{3}+\pi\eta^3\right)\sin^2{\beta}
    +\left(\frac{4\sin{\gamma}~\eta^2}{3}+\frac{\pi\cos^2{\beta}~\eta^3}{2}+\frac{4\pi\eta^4}{15}
    \right)\right\}} \\ [0.3cm]
    \displaystyle{L^2D^3\left(\frac{2\sin{\gamma}}{3}+\frac{\pi\eta}{2}+\frac{4\pi\eta^2}{15}\right)}
\end{array}\right].
\]

Following the similar process, we can also work out the fourth moment, which is used to modeling the smectic phase.
In region A, we have
\begin{align*}
\displaystyle{\int_\text{A} G(|\ur|,\um,\um')
\left[\begin{array}{c}
    r_1^4 \\ [0.06cm]
    r_2^4\\ [0.06cm]
    r_3^4
\end{array}\right]\ud\ur}
=~&\sin{\gamma}\int_{-D}^{D}\int_{-L/2}^{L/2}\int_{-L/2}^{L/2}
\left[\begin{array}{c}
    (u+v)^4\cos^4{\beta}  \\ [0.06cm]
    (u-v)^4\sin^4{\beta}\\ [0.06cm]
    z^4
\end{array}\right]~\ud u\ud v\ud z \\
=~&\left[\begin{array}{c}
    2L^6D\cos^4{\beta}\sin{\gamma}/15 \\ [0.06cm]
    2L^6D\sin^4{\beta}\sin{\gamma}/15 \\ [0.06cm]
    2L^2D^5\sin{\gamma}/5
\end{array}\right],
\end{align*}

\begin{align*}
\displaystyle{\int_\text{A} G(|\ur|,\um,\um')
\left[\begin{array}{c}
    r_1^2r_2^2 \\ [0.06cm]
    r_2^2r_3^2  \\ [0.06cm]
    r_1^2r_3^2
\end{array}\right]\ud\ur}
=~&\sin{\gamma}\int_{-D}^{D}\int_{-L/2}^{L/2}\int_{-L/2}^{L/2}
\left[\begin{array}{c}
    (u^2-v^2)^2\cos^2{\beta}\sin^2{\beta}  \\ [0.06cm]
    (u-v)^2\sin^2{\beta}z^2  \\ [0.06cm]
    (u+v)^2\cos^2{\beta}z^2
\end{array}\right]~\ud u\ud v\ud z \\
=~&\left[\begin{array}{c}
    L^6D\cos^2{\beta}\sin^2{\beta}\sin{\gamma}/45 \\ [0.06cm]
    L^4D^3\sin^2{\beta}\sin{\gamma}/9 \\ [0.06cm]
    L^4D^3\cos^2{\beta}\sin{\gamma}/9
\end{array}\right].
\end{align*}

In region B, we get
\begin{align*}
&\displaystyle{\int_\text{B} G(\um,\um',\ur)
\left[\begin{array}{c}
    r_1^4 \\ [0.06cm]
    r_2^4\\ [0.06cm]
    r_3^4
\end{array}\right]\ud\ur} \\
=~&4\int_{0}^{L}\int_{0}^{D}\int_{0}^{\pi}
r\left[\begin{array}{c}
    (t\cos{\beta}+r\sin{\beta}\sin{\theta})^4  \\ [0.06cm]
    ((L-t)\sin{\beta}+r\cos{\beta}\sin{\theta})^4 \\ [0.06cm]
    (r\cos{\theta})^4  \\ [0.16cm]
\end{array}\right]~\ud \theta\ud r\ud t \\
=~&\left[\begin{array}{c}
    2\pi L^5D^2\cos^4{\beta}/5+8L^4D^3\cos^3{\beta}\sin{\beta}/3+\pi L^3D^4\cos^2{\beta}\sin^2{\beta}  \\ [0.06cm]
    +32L^2D^5\cos{\beta}\sin^3{\beta}/15+\pi LD^6\sin^4{\beta}/4 \\ [0.26cm]
    2\pi L^5D^2\sin^4{\beta}/5+8L^4D^3\sin^3{\beta}\cos{\beta}/3+\pi L^3D^4\sin^2{\beta}\cos^2{\beta}  \\ [0.06cm]
    +32L^2D^5\sin{\beta}\cos^3{\beta}/15+\pi LD^6\cos^4{\beta}/4 \\ [0.26cm]
    \pi LD^6/4
\end{array}\right],
\end{align*}

\begin{align*}
&\displaystyle{\int_\text{B} G(\um,\um',\ur)
\left[\begin{array}{c}
    r_1^2r_2^2 \\ [0.06cm]
    r_2^2r_3^2\\ [0.06cm]
    r_1^2r_3^2
\end{array}\right]\ud\ur} \\
=~&4\int_{0}^{L}\int_{0}^{D}\int_{0}^{\pi}
r\left[\begin{array}{c}
    (t\cos{\beta}+r\sin{\beta}\sin{\theta})^2((L-t)\sin{\beta}+r\cos{\beta}\sin{\theta})^2  \\ [0.06cm]
    ((L-t)\sin{\beta}+r\cos{\beta}\sin{\theta})^2(r\cos{\theta})^2 \\ [0.06cm]
    (t\cos{\beta}+r\sin{\beta}\sin{\theta})^2(r\cos{\theta})^2 \\ [0.16cm]
\end{array}\right]~\ud \theta\ud r\ud t \\
=~&\left[\begin{array}{c}
    \pi L^3D^4(\cos^4{\beta}+\sin^4{\beta})/6+(4L^4D^3/9+16L^2D^5/15)(\cos^3{\beta}\sin{\beta}+\sin^3{\beta}\cos{\beta})  \\ [0.06cm]
    +(\pi L^5D^2/15+\pi L^3D^4/3+\pi LD^6/4)\cos^2{\beta}\sin^2{\beta} \\ [0.26cm]
    \pi L^3D^4\sin^2{\beta}/6+8 L^2D^5\cos{\beta}\sin{\beta}/15 +\pi LD^6\cos^2{\beta}/12 \\ [0.26cm]
    \pi L^3D^4\cos^2{\beta}/6+8 L^2D^5\cos{\beta}\sin{\beta}/15 +\pi LD^6\sin^2{\beta}/12 \\
\end{array}\right].
\end{align*}

Finally in region C\Rmnum{1} and C\Rmnum{2}, it reads
\begin{align*}
&\displaystyle{\int_\text{CI} G(|\ur|,\um,\um')
\left[\begin{array}{c}
    r_1^4 \\ [0.06cm]
    r_2^4\\ [0.06cm]
    r_3^4
\end{array}\right]\ud\ur} \\
=~&\int_{0}^{D}\int_{0}^{\pi}\int_{\pi/2-\beta}^{\pi/2+\beta}
r^2\sin{\theta}\left[\begin{array}{c}
    (r\sin{\theta}\cos{\varphi})^4  \\ [0.06cm]
    (L\sin{\beta}+r\sin{\theta}\sin{\varphi})^4 \\ [0.06cm]
    (r\cos{\theta})^4
\end{array}\right]~\ud \varphi\ud \theta\ud r \\
=~&\left[\begin{array}{c}
    16D^7(3\beta/4+\cos^3{\beta}\sin{\beta}/4-\sin^3{\beta}\cos{\beta}/4-\cos{\beta}\sin{\beta})/105  \\ [0.26cm]
    2L^4D^3(2\beta)\sin^4{\beta}/3+\pi L^3D^4\sin^4{\beta}+8L^2D^5\sin^2{\beta}(\beta+\cos{\beta}\sin{\beta})/5  \\ [0.03cm]
    +\pi LD^6\sin^2{\beta}(3-\sin^2{\beta})/6  \\ [0.03cm]
    +16D^7(3\beta/4+\cos^3{\beta}\sin{\beta}/4-\sin^3{\beta}\cos{\beta}/4+\cos{\beta}\sin{\beta})/105 \\ [0.26cm]
    4D^7\beta/35
\end{array}\right],
\end{align*}

\begin{align*}
&\displaystyle{\int_\text{CI} G(|\ur|,\um,\um')
\left[\begin{array}{c}
    r_1^2r_2^2 \\ [0.06cm]
    r_2^2r_3^2\\ [0.06cm]
    r_1^2r_3^2
\end{array}\right]\ud\ur} \\
=~&\int_{0}^{D}\int_{0}^{\pi}\int_{\pi/2-\beta}^{\pi/2+\beta}
r^2\sin{\theta}\left[\begin{array}{c}
    (r\sin{\theta}\cos{\varphi})^2(L\sin{\beta}+r\sin{\theta}\sin{\varphi})^2  \\ [0.06cm]
    (L\sin{\beta}+r\sin{\theta}\sin{\varphi})^2(r\cos{\theta})^2 \\ [0.06cm]
    (r\cos{\theta})^2(r\sin{\theta}\cos{\varphi})^2
\end{array}\right]~\ud \varphi\ud \theta\ud r \\
=~&\left[\begin{array}{c}
    4L^2D^5\sin^2{\beta}(\beta-\cos{\beta}\sin{\beta})/15+\pi LD^6\sin^4{\beta}/12 \\ [0.03cm]
    +16D^7(\beta/4-\cos^3{\beta}\sin{\beta}/4+\sin^3{\beta}\cos{\beta}/4)/105  \\ [0.26cm]
    2L^2D^5(2\beta)\sin^2{\beta}/15+\pi LD^6\sin^2{\beta}/12+4D^7(\beta+\cos{\beta}\sin{\beta})/105  \\ [0.26cm]
    4D^7(\beta-\cos{\beta}\sin{\beta})/105
\end{array}\right],
\end{align*}
\begin{align*}
&\displaystyle{\int_\text{CII} G(|\ur|,\um,\um')
\left[\begin{array}{c}
    r_1^4 \\ [0.06cm]
    r_2^4\\ [0.06cm]
    r_3^4
\end{array}\right]\ud\ur} \\
=~&\int_{0}^{D}\int_{0}^{\pi}\int_{-\pi/2+\beta}^{\pi/2-\beta}
r^2\sin{\theta}\left[\begin{array}{c}
    (L\cos{\beta}+r\sin{\theta}\cos{\varphi})^4  \\ [0.06cm]
    (r\sin{\theta}\sin{\varphi})^4\\ [0.06cm]
    (r\cos{\theta})^4  \\ [0.16cm]
\end{array}\right]~\ud \varphi\ud \theta\ud r \\
=~&\left[\begin{array}{c}
    2L^4D^3(\pi-2\beta)\cos^4{\beta}/3+\pi L^3D^4\cos^4{\beta}+8L^2D^5\cos^2{\beta}(\pi/2-\beta+\cos{\beta}\sin{\beta})/5  \\ [0.03cm]
    +\pi LD^6\cos^2{\beta}(3-\cos^2{\beta})/6  \\ [0.03cm]
    +16D^7(3\pi/8-3\beta/4-\cos^3{\beta}\sin{\beta}/4+\sin^3{\beta}\cos{\beta}/4+\cos{\beta}\sin{\beta})/105 \\ [0.26cm]
    16D^7(3\pi/8-3\beta/4-\cos^3{\beta}\sin{\beta}/4+\sin^3{\beta}\cos{\beta}/4-\cos{\beta}\sin{\beta})/105  \\ [0.26cm]
    2D^7(\pi-2\beta)/35
\end{array}\right],
\end{align*}

\begin{align*}
&\displaystyle{\int_\text{CII} G(|\ur|,\um,\um')
\left[\begin{array}{c}
    r_1^2r_2^2 \\ [0.06cm]
    r_2^2r_3^2\\ [0.06cm]
    r_1^2r_3^2
\end{array}\right]\ud\ur} \\
=~&\int_{0}^{D}\int_{0}^{\pi}\int_{-\pi/2+\beta}^{\pi/2-\beta}
r^2\sin{\theta}\left[\begin{array}{c}
    (L\cos{\beta}+r\sin{\theta}\cos{\varphi})^2(r\sin{\theta}\sin{\varphi})^2  \\ [0.06cm]
    (r\sin{\theta}\sin{\varphi})^2(r\cos{\theta})^2\\ [0.06cm]
    (r\cos{\theta})^2(L\cos{\beta}+r\sin{\theta}\cos{\varphi})^2  \\ [0.16cm]
\end{array}\right]~\ud \varphi\ud \theta\ud r \\
=~&\left[\begin{array}{c}
    4L^2D^5\cos^2{\beta}(\pi/2-\beta-\cos{\beta}\sin{\beta})/15+\pi LD^6\cos^4{\beta}/12 \\ [0.03cm]
    +16D^7(\pi/8-\beta/4+\cos^3{\beta}\sin{\beta}/4-\sin^3{\beta}\cos{\beta}/4)/105  \\ [0.26cm]
    4D^7(\pi/2-\beta-\cos{\beta}\sin{\beta})/105  \\ [0.26cm]
    2L^2D^5(\pi-2\beta)\cos^2{\beta}/15+\pi LD^6\cos^2{\beta}/12+4D^7(\pi/2-\beta+\cos{\beta}\sin{\beta})/105
\end{array}\right].
\end{align*}

Sum them up, and finally we gain the entire fourth moment result as follows:
\begin{align*}
&\displaystyle{\int G(|\ur|,\um,\um') \left[\begin{array}{c}
    r_1^4 \\ [0.06cm]
    r_2^4\\ [0.06cm]
    r_3^4
\end{array}\right]~\ud\ur} \\
=~&\left[\begin{array}{c}
     \displaystyle{L^6D\bigg\{\bigg(\frac{2\sin{\gamma}}{15}+\frac{2\pi \eta}{5}+\frac{4(\pi-\gamma)\eta^2}{3}+2\pi\eta^3-\frac{\pi \eta^5}{12}\bigg)\cos^4{\beta}+\bigg(\frac{8\eta^2}{3}+\frac{16\eta^4}{15}\bigg)\cos^3{\beta}\sin{\beta}}\\ [0.03cm]
     \displaystyle{\bigg(\pi\eta^3+\frac{8(\pi-\gamma)\eta^4}{5}+\frac{\pi\eta^5}{2}\bigg)\cos^2{\beta}+\frac{32\eta^4}{15}\cos{\beta}\sin{\beta}
     +\bigg(\frac{\pi\eta^5}{4}+\frac{4\pi\eta^6}{35}\bigg)\bigg\}}\\ [0.6cm]
     \displaystyle{L^6D\bigg\{\bigg(\frac{2\sin{\gamma}}{15}+\frac{2\pi \eta}{5}+\frac{4\gamma\eta^2}{3}+2\pi\eta^3-\frac{\pi \eta^5}{12}\bigg)\sin^4{\beta}+\bigg(\frac{8\eta^2}{3}+\frac{16\eta^4}{15}\bigg)\sin^3{\beta}\cos{\beta}}\\ [0.03cm]
     \displaystyle{\bigg(\pi\eta^3+\frac{8\gamma\eta^4}{5}+\frac{\pi\eta^5}{2}\bigg)\sin^2{\beta}+\frac{32\eta^4}{15}\sin{\beta}\cos{\beta}
     +\bigg(\frac{\pi\eta^5}{4}+\frac{4\pi\eta^6}{35}\bigg)\bigg\}}\\ [0.6cm]
     \displaystyle{L^6D\bigg\{\frac{2\sin{\gamma}~\eta^4}{5}+\frac{\pi\eta^5}{4}+\frac{4\pi\eta^6}{35}\bigg\}}\\
\end{array}\right],
\end{align*}

\begin{align*}
&\displaystyle{\int G(|\ur|,\um,\um') \left[\begin{array}{c}
    r_1^2r_2^2 \\ [0.06cm]
    r_2^2r_3^2\\ [0.06cm]
    r_1^2r_3^2
\end{array}\right]~\ud\ur} \\
=~&\left[\begin{array}{c}
     \displaystyle{L^6D\bigg\{\bigg(\frac{\sin{\gamma}}{45}+\frac{\pi\eta}{15}-\frac{\pi\eta^5}{12}\bigg)\cos^2{\beta}\sin^2{\beta}
     +\bigg(\frac{4\eta^2}{9}+\frac{8\eta^4}{15}\bigg)\cos{\beta}\sin{\beta}}\\ [0.03cm]
     \displaystyle{+\frac{4\eta^4}{15}(\gamma\sin^2{\beta}+(\pi-\gamma)\cos^2{\beta})+\bigg(\frac{\pi\eta^3}{6}+\frac{\pi\eta^5}{6}+\frac{4\pi\eta^6}{105}\bigg)
     \bigg\}}\\ [0.6cm]
     \displaystyle{L^6D\bigg\{\bigg(\frac{\sin{\gamma}~\eta^2}{9}+\frac{\pi\eta^3}{6}+\frac{4\gamma\eta^4}{15}
     +\frac{\pi\eta^5}{12}\bigg)\sin^2{\beta}
     +\frac{8\eta^4}{15}\cos{\beta}\sin{\beta}+\bigg(\frac{\pi\eta^5}{12}+\frac{4\pi\eta^6}{105}\bigg)\bigg\}}\\ [0.6cm]
     \displaystyle{L^6D\bigg\{\bigg(\frac{\sin{\gamma}~\eta^2}{9}+\frac{\pi\eta^3}{6}+\frac{4(\pi-\gamma)\eta^4}{15}
     +\frac{\pi\eta^5}{12}\bigg)
     \cos^2{\beta}+\frac{8\eta^4}{15}\cos{\beta}\sin{\beta}+\bigg(\frac{\pi\eta^5}{12}+\frac{4\pi\eta^6}{105}\bigg)\bigg\}}\\ [0.6cm]
\end{array}\right].
\end{align*}

Similar to the decomposition of second moment, we let
\begin{equation*}
\begin{array}{l}
\displaystyle{~~~M^{(4)}} \\ [0.26cm]
\displaystyle{=W_1(\um,\um')~\nn_1\nn_1\nn_1\nn_1+W_2(\um,\um')~\nn_2\nn_2\nn_2\nn_2+W_3(\um,\um')~\nn_3\nn_3\nn_3\nn_3}  \\ [0.1cm]
\displaystyle{\quad\quad+W_4(\um,\um')~(\nn_1\nn_1\nn_2\nn_2)_{\mathrm{sym}}+W_5(\um,\um')~(\nn_2\nn_2\nn_3\nn_3)_{\mathrm{sym}}} \\ [0.1cm]
\displaystyle{\quad\quad+W_6(\um,\um')~(\nn_1\nn_1\nn_3\nn_3)_{\mathrm{sym}}}  \\ [0.1cm]
\displaystyle{=\frac{W_3}{3}(\udeltayier\udeltasansi)_{\mathrm{sym}}+(W_1+W_3-6W_6)~\nn_1\nn_1\nn_1\nn_1} \\ [0.1cm]
\displaystyle{\quad\quad+(W_2+W_3-6W_5)~\nn_2\nn_2\nn_2\nn_2+(W_4+\frac{W_3}{3}-W_5-W_6)(\nn_1\nn_1\nn_2\nn_2)_{\mathrm{sym}}}  \\ [0.1cm]
\displaystyle{\quad\quad+(W_6-\frac{W_3}{3})(\udeltayier\nn_1\nn_1)_{\mathrm{sym}}+
(W_5-\frac{W_3}{3})(\udeltayier\nn_2\nn_2)_{\mathrm{sym}}.}  \\ [0.3cm]
\displaystyle{=\frac{W_3}{3}(\udeltayier\udeltasansi)_{\mathrm{sym}}+\frac{W_1+W_3-6W_6}{16\cos^4{\beta}}\bigg[(\um\um\um\um+\um'\um'\um'\um')
+(\um\um\um'\um')_{\mathrm{sym}}}  \\ [0.30cm]
\displaystyle{\quad\quad+(\um\um\um\um'+\um'\um'\um'\um)_{\mathrm{sym}}\bigg]+\frac{W_2+W_3-6W_5}{16\sin^4{\beta}}\bigg[(\um\um\um\um+\um'\um'\um'\um')} \\ [0.30cm]
\displaystyle{\quad\quad
+(\um\um\um'\um')_{\mathrm{sym}}-(\um\um\um\um'+\um'\um'\um'\um)_{\mathrm{sym}}\bigg]}  \\ [0.30cm]
\displaystyle{\quad\quad+\frac{W_3/3+W_4-W_5-W_6}{16\cos^2{\beta}\sin^2{\beta}}\bigg[6(\um\um\um\um+\um'\um'\um'\um')-2(\um\um\um'\um')_{\mathrm{sym}}\bigg]} \\ [0.30cm]
\displaystyle{\quad\quad+\frac{W_6-W_3/3}{4\cos^2{\beta}}\bigg[(\udeltayier\um\um+\udeltayier\um'\um')_{\mathrm{sym}}+(\udeltayier\um\um')_{\mathrm{sym}}\bigg]} \\ [0.30cm]
\displaystyle{\quad\quad+\frac{W_5-W_3/3}{4\sin^2{\beta}}\bigg[(\udeltayier\um\um+\udeltayier\um'\um')_{\mathrm{sym}}-(\udeltayier\um\um')_{\mathrm{sym}}\bigg]} \\ [0.30cm]
\end{array} \\
\end{equation*}
Direct simplification gives us (\ref{fourthmoment}).

\subsection{The derivation of the elastic energy in the $Q$-tensor form}
Dropping the high order terms, we can get
\begin{align}
\frac{4}{\pi L^4Dk_BT}&F^{(2)}_{\text{elastic}}=-\int\bigg\{
\Big(\alpha_{11}+\alpha_{12}P_2(\mm\cdot\mm')+\alpha_{13}P_4(\mm\cdot\mm')\Big)
\partial_i f(\xx,\mm')\partial_i f(\xx,\mm)
\nonumber\\
&+\Big(\alpha_{21}-\frac12\alpha_{22}+\alpha_{22}\frac{3}{2}(\mm\cdot\mm')^2\Big)
(\mm\mm+\mm'\mm'):\nabla  f(\xx,\mm')\nabla  f(\xx,\mm)\nonumber\\
&+\Big(\big(\alpha_{31}-\frac{1}{2}\alpha_{32}\big)\mm\cdot\mm'+\frac32\alpha_{32}(\mm\cdot\mm')^2\Big)
(\mm\mm'+\mm'\mm):\nabla  f(\xx,\mm')\nabla  f(\xx,\mm)\bigg\}\ud\mm'\ud\mm\nonumber\\
=&-\int\bigg\{
\Big(\alpha_{11}+\frac{3}{2}\alpha_{12}\Xi_2:\Xi_2'+\frac{35}{8}\alpha_{13}\Xi_4:\Xi_4'\Big)
\partial_i f(\xx,\mm')\partial_i f(\xx,\mm)
\nonumber\\
&+\Big(\alpha_{21}(\Xi_2+\Xi_2'+\frac23\II)+3\alpha_{22}\mm\mm\mm\mm:\Xi_2'\Big)\nabla  f(\xx,\mm')\nabla  f(\xx,\mm)\nonumber\\
&+\big(\alpha_{31}-\frac{1}{2}\alpha_{32}\big)(\Xi_2\cdot \Xi_2'
+\Xi_2'\cdot \Xi_2+\frac{2}{3}\Xi_2+\frac23\Xi_2'+\frac29\II):\nabla  f(\xx,\mm')\nabla  f(\xx,\mm)\nonumber\\
&+\frac32\alpha_{32}(\mm\cdot\mm')^3
(\mm\mm'+\mm'\mm):\nabla  f(\xx,\mm')\nabla  f(\xx,\mm)\bigg\}\ud\mm'\ud\mm
\end{align}
It is direct to check that
\begin{align*}
&(\mm\cdot\mm')\mm\mm'=\Xi_2\cdot \Xi_2'+\frac{1}{3}(\Xi_2+\Xi_2')+\frac19\II.\\
&\mm\mm\mm\mm:\Xi_2'=\Xi_4:\Xi_2'+\frac{1}{7}\Big(2\Xi_2\cdot \Xi_2'+2\Xi_2'\cdot \Xi_2+\II \Xi_2:\Xi_2'\Big)+\frac{2}{15}\Xi_2',\\
&(\mm\cdot\mm')^2(\mm\mm+\mm'\mm')=\mm\mm\mm\mm:\Xi_2'+\mm'\mm'\mm'\mm':\Xi_2+\frac13(\Xi_2+\Xi_2'+\frac23\II),
\end{align*}
and
\begin{align}
&\Xi_{4,\alpha\gamma\mu\nu}\Xi_{4,\beta\gamma\mu\nu}'\nonumber\\
&=(\mm\cdot\mm')^3\mm\mm'-\frac{15}{49}(\mm\cdot\mm')\mm\mm'
+\frac{6}{49}(\mm\cdot\mm')\mm'\mm-\frac37(\mm\cdot\mm')^2(\mm\mm+\mm'\mm')\nonumber\\
&\qquad+\frac{3}{49}(\mm\mm+\mm'\mm')
+\frac{3}{49}(\mm\cdot\mm')^2\II-\frac{3\II}{7\cdot35}\nonumber\\
&=(\mm\cdot\mm')^3\mm\mm'-\frac{15}{49}\Xi_2\cdot \Xi_2'
+\frac{6}{49}\Xi_2'\cdot \Xi_2-\frac37(\mm\cdot\mm')^2(\mm\mm+\mm'\mm')+\frac{\II}{49}(3\Xi_2:\Xi_2'+\frac75).\nonumber
\end{align}
Hence we have
\begin{align}
&\frac{4}{\pi L^4Dk_BT}F^{(2)}_{\text{elastic}}\nonumber\\
=&-\int\bigg\{
\Big(\alpha_{11}+\frac{3}{2}\alpha_{12}\Xi_2:\Xi_2'+\frac{35}{8}\alpha_{13}\Xi_4:\Xi_4'\Big)
\partial_i f(\xx,\mm')\partial_i f(\xx,\mm)
\nonumber\\
&+\Big(\alpha_{21}(\Xi_2+\Xi_2'+\frac23\II)+3\alpha_{22}\mm\mm\mm\mm:\Xi_2'\Big)\nabla  f(\xx,\mm')\nabla  f(\xx,\mm)\nonumber\\
&+\big(\alpha_{31}-\frac{1}{2}\alpha_{32}\big)\Big(\Xi_2\cdot \Xi_2'
+\Xi_2'\cdot \Xi_2+\frac{2}{3}\Xi_2+\frac23\Xi_2'+\frac29\II\Big):\nabla  f(\xx,\mm')\nabla  f(\xx,\mm)\nonumber\\
&+\frac32\alpha_{32}\Big(\Xi_4\dot{:}\Xi_4'+\Xi_4'\dot{:}\Xi_4+\frac{9}{49}\Xi_2\cdot \Xi_2'
+\frac{9}{49}\Xi_2'\cdot \Xi_2+\frac67(\mm\cdot\mm')^2(\mm\mm+\mm'\mm')\nonumber\\
&\qquad
-\frac{2\II}{49}(3\Xi_2:\Xi_2'+\frac75)\Big):\nabla  f(\xx,\mm')\nabla  f(\xx,\mm)\bigg\}\ud\mm'\ud\mm\nonumber\\
=&-\int\bigg\{
\Big(\alpha_{11}+\frac{3}{2}\alpha_{12}\Xi_2:\Xi_2'+\frac{35}{8}\alpha_{13}\Xi_4:\Xi_4'\Big)
\partial_i f(\xx,\mm')\partial_i f(\xx,\mm)
\nonumber\\
&+\Big(\alpha_{21}(\Xi_2+\Xi_2'+\frac23\II)\Big)\nabla  f(\xx,\mm')\nabla  f(\xx,\mm)\nonumber\\
&+\big(\alpha_{31}-\frac{1}{2}\alpha_{32}\big)\Big(\Xi_2\cdot \Xi_2'
+\Xi_2'\cdot \Xi_2+\frac{2}{3}\Xi_2+\frac23\Xi_2'+\frac29\II\Big):\nabla  f(\xx,\mm')\nabla  f(\xx,\mm)\nonumber\\
&+\frac32\alpha_{32}\Big(\Xi_4\dot{:}\Xi_4'+\Xi_4'\dot{:}\Xi_4+\frac{9}{49}\Xi_2\cdot \Xi_2'
+\frac{9}{49}\Xi_2'\cdot \Xi_2-\frac{2\II}{49}(3\Xi_2:\Xi_2'+\frac75)+\frac27(\Xi_2+\Xi_2'+\frac23\II)\Big):\nabla  f\nabla  f'\nonumber\\
&+\big(3\alpha_{22}+\frac{9}{7}\alpha_{32}\big)\Big(\Xi_4:\Xi_2'
+\frac{1}{7}\big(2\Xi_2\cdot \Xi_2'+2\Xi_2'\cdot \Xi_2+\II \Xi_2:\Xi_2'\big)+\frac{2}{15}
\Xi_2'\Big)\bigg\}:\nabla  f\nabla  f'\ud\mm'\ud\mm\nonumber\\
=&-\bigg\{\Big(\alpha_{11}+\frac23\alpha_{21}+\frac29\alpha_{31}-\frac{2}{18}
\alpha_{32}+\frac1{5}\alpha_{32}\Big)|\nabla{c}|^2\nonumber\\
&+\Big(\frac32\alpha_{12}-\frac9{49}\alpha_{32}+\frac37\alpha_{22}+\frac{18}{49}\alpha_{32}\Big)
|\nabla(cQ_2)|^2+\frac{35}{8}\alpha_{13}|\nabla(cQ_4)|^2\nonumber\\
&+\Big(2\alpha_{21}+\frac{4}{3}(\alpha_{31}-\frac{1}{2}\alpha_{32})+\frac67\alpha_{32}
+\frac{2}{15}(3\alpha_{22}+\frac{18}{7}\alpha_{32})\Big)\partial_i(cQ_{2ij})\partial_jc\nonumber\\
&+\Big(\alpha_{31}-\frac{1}{2}\alpha_{32}+\frac{27}{98}\alpha_{32}
+\frac27(3\alpha_{22}+\frac{18}{7}\alpha_{32})\Big)
\Big(\partial_i(cQ_{ik})\partial_j(cQ_{jk})+\partial_i(cQ_{jk})\partial_j(cQ_{ik})\Big)\nonumber\\
&+\frac32\alpha_{32}\Big(\partial_i(cQ_{4iklm})\partial_j
(cQ_{4jklm})+\partial_i(cQ_{4jklm})\partial_j(cQ_{4iklm})\Big)\nonumber\\
&+\big(3\alpha_{22}+\frac{18}{7}\alpha_{32}\big)\partial_i(cQ_{4ijkl})\partial_j(cQ_{2kl})\bigg\}.\label{Q-Energy0}
\end{align}
Finally, we get
\begin{align}
&\frac{4}{\pi L^4Dk_BT}F^{(2)}_{\text{elastic}}\nonumber\\
=&-\bigg\{\Big(\alpha_{11}+\frac23\alpha_{21}+\frac29\alpha_{31}+\frac{4}{45}
\alpha_{32}\Big)|\nabla{c}|^2+\Big(\frac32\alpha_{12}+\frac37\alpha_{22}
+\frac{9}{49}\alpha_{32}\Big)|\nabla(cQ_2)|^2\nonumber\\
&+\frac{35}{8}\alpha_{13}|\nabla(cQ_4)|^2+\Big(2\alpha_{21}+\frac{4}{3}
\alpha_{31}+\frac{2}{5}\alpha_{22}+\frac8{15}\alpha_{32}
\Big)\partial_i(cQ_{2ij})\partial_jc\nonumber\\
&+\Big(\alpha_{31}+\frac{25}{49}\alpha_{32}+\frac67\alpha_{22}\Big)
\Big(\partial_i(cQ_{ik})\partial_j(cQ_{jk})+\partial_i(cQ_{jk})\partial_j(cQ_{ik})\Big)\nonumber\\
&+\frac32\alpha_{32}\Big(\partial_i(cQ_{4iklm})\partial_j
(cQ_{4jklm})+\partial_i(cQ_{4jklm})\partial_j(cQ_{4iklm})\Big)\nonumber\\
&+\big(3\alpha_{22}+\frac{18}{7}\alpha_{32}\big)\partial_i(cQ_{4ijkl})\partial_j(cQ_{2kl})\bigg\}.
\end{align}
The calculation of fourth moment gives us
\begin{align*}
&\frac{24}{\pi L^6D k_BT}F^{(4)}_{\text{elastic}}\nonumber\\
\approx &\int_{\Omega}\int_{\BS}\int_{\BS} \bigg\{ 2\mu_{11}~f(\xx,\mm)~m_im_jm_km_l~\partial_{ijkl}\big\{f(\xx,\mm')\big\}\ud\mm'\ud\mm  \nonumber\\
&\quad\quad\quad\quad\quad+(\mu_{21}-\frac{1}{2}\mu_{22})f(\xx,\mm)(m_im_jm'_km'_l)_{\mathrm{sym}}~\partial_{ijkl}\big\{f(\xx,\mm')\big\}
\ud\mm'\ud\mm \nonumber\\
&\quad\quad\quad\quad\quad+\frac{3}{2}\mu_{22}f(\xx,\mm)(\mm\cdot\mm')^2(m_im_jm'_km'_l)_{\mathrm{sym}}~\partial_{ijkl}\big\{f(\xx,\mm')\big\}
\ud\mm'\ud\mm\bigg\} \ud\xx.
\end{align*}

With periodic boundary condition, we get
\begin{align*}
&\int c(\xx)\rho(\xx,\mm)~m_im_jm_km_l~\partial_{ijkl}(c(\xx)\rho(\xx,\mm'))~\ud\mm'\ud\mm\ud\xx \nonumber \\
=&\int c~\partial_{ijkl}\big\{c(Q_{4ijkl}+\frac{1}{7}(Q_{2ij}\delta_{kl})_{\mathrm{symmtric}}+\frac{1}{15}(\delta_{ij}\delta_{kl})_{\mathrm{symmtric}})\big\}\ud\xx \nonumber \\
=&\int \partial_{ij}(cQ_{4ijkl})\partial_{kl}(c)+\frac{1}{7}c(\partial_{ijkk}(cQ_{2ij}))_{\mathrm{symmtric}}+\frac{1}{15}c(\partial_{iikk}(c))_{\mathrm{symmtric}}
~\ud\xx \nonumber \\
=&\int \bigg\{\partial_{ij}(cQ_{4ijkl})\partial_{kl}(c)+\frac{6}{7}\partial_{ij}(cQ_{2ij})\partial_{kk}(c)
+\frac{1}{5}(\partial_{ij}(c))^2\bigg\}\ud\xx,
\end{align*}
and
\begin{align*}
&\int c(\xx)\rho(\xx,\mm)(m_im_jm'_km'_l)_{\mathrm{symmtric}}\partial_{ijkl}(c(\xx)\rho(\xx,\mm'))~\ud\mm'\ud\mm\ud\xx \nonumber \\
=&\int \big[c(Q_{2ij}+\frac{1}{3}\delta_{ij})~\partial_{ijkl}(c(Q_{2kl}+\frac{1}{3}\delta_{kl}))\big]_{\mathrm{symmtric}}~\ud\xx \nonumber \\
=&\int \bigg\{6\partial_{ij}(cQ_{2ij})\partial_{kl}(cQ_{2kl})+4\partial_{ij}(cQ_{2ij})\partial_{kk}(c)+\frac{2}{3}(\partial_{ij}(c))^2\bigg\}\ud\xx.
\end{align*}
In addition, we can derive that
\begin{align*}
&\int c(\xx)\rho(\xx,\mm)(\mm\cdot\mm')^2(m_im_jm'_km'_l)_{\mathrm{symmtric}}\partial_{ijkl}(c(\xx)\rho(\xx,\mm'))~\ud\mm'\ud\mm\ud\xx \nonumber \\
=&6\int c\Big(Q_{4ijkl}+\frac{1}{7}(Q_{2ij}\delta_{kl})_{\mathrm{symmtric}}
+\frac{1}{15}(\delta_{ij}\delta_{kl})_{\mathrm{symmtric}}\Big) \nonumber \\
&~~~~\cdot\partial_{ijkl}\Big\{(c(Q_{4ijkl}+\frac{1}{7}(Q_{2ij}\delta_{kl})_{\mathrm{symmtric}}
+\frac{1}{15}(\delta_{ij}\delta_{kl})_{\mathrm{symmtric}})\Big\}~\ud\xx \nonumber \\
=&\int \bigg\{6\partial_{ij}(cQ_{4ijpq})\partial_{kl}(cQ_{4klpq})
+\Big(\frac{114}{49}\partial_{ij}(cQ_{2ij})\partial_{kl}(cQ_{2kl})+\frac{96}{49}\partial_{ik}(cQ_{2ip})\partial_{jk}(cQ_{2jp}) \nonumber \\
&~~~~+\frac{12}{49}\partial_{ij}(cQ_{2ij})\partial_{kk}(cQ_{2pp})+\frac{6}{49}\partial_{ij}(cQ_{2pq})\partial_{ij}(cQ_{2pq})\Big)
+\frac{22}{75}\partial_{ij}(c)\partial_{ij}(c) \nonumber \\
&~~~~+\Big(\frac{12}{7}\partial_{ij}(cQ_{4ijpp})\partial_{kl}(cQ_{2kl})+\frac{48}{7}\partial_{ij}(cQ_{4ijkp})\partial_{kl}(cQ_{2lp})
+\frac{12}{7}\partial_{ij}(cQ_{4ijpq})\partial_{kk}(cQ_{2pq})\Big) \nonumber \\
&~~~~+\Big(\frac{4}{5}\partial_{ij}(cQ_{4ijpp})\partial_{kk}(c)+\frac{8}{5}\partial_{ij}(cQ_{4ijkl})\partial_{kl}(c)\Big) \nonumber \\
&~~~~+\Big(\frac{76}{35}\partial_{ij}(cQ_{2ij})\partial_{kk}(c)+\frac{4}{35}\partial_{ij}(cQ_{2pp})\partial_{ij}(c)\Big)\bigg\}\ud\xx. \nonumber
\end{align*}
Together with the trace free property of Q-tensor, we sum up to obtain the fourth term as
\begin{align}
&\frac{24}{\pi L^6D k_BT}F^{(4)}_{\text{elastic}}\nonumber\\
\approx &\int \bigg\{\big(9\mu_{22}\big)\partial_{ij}(cQ_{4ijpq})\partial_{kl}(cQ_{4klpq})
+\Big(\big(6\mu_{21}+\frac{24}{49}\mu_{22}\big)\partial_{ij}(cQ_{2ij})\partial_{kl}(cQ_{2kl}) \nonumber \\
&~~~~+\big(\frac{144}{49}\mu_{22}\big)\partial_{ik}(cQ_{2ip})\partial_{jk}(cQ_{2jp})
+\big(\frac{9}{49}\mu_{22}\big)\partial_{ij}(cQ_{2pq})\partial_{ij}(cQ_{2pq})\Big) \nonumber \\
&~~~~+\big(\frac{2}{5}\mu_{11}+\frac{2}{3}\mu_{21}+\frac{8}{75}\mu_{22}\big)\partial_{ij}(c)\partial_{ij}(c)
+\Big(\big(\frac{72}{7}\mu_{22}\big)\partial_{ij}(cQ_{4ijkp})\partial_{kl}(cQ_{2lp}) \nonumber \\
&~~~~+\big(\frac{18}{7}\mu_{22}\big)\partial_{ij}(cQ_{4ijpq})\partial_{kk}(cQ_{2pq})\Big)
+\big(2\mu_{11}+\frac{12}{5}\mu_{22}\big)\partial_{ij}(cQ_{4ijkl})\partial_{kl}(c)  \nonumber \\
&~~~~+\big(\frac{12}{7}\mu_{11}+4\mu_{21}+\frac{44}{35}\mu_{22}\big)\partial_{ij}(cQ_{2ij})\partial_{kk}(c)\bigg\}\ud\xx.  \label{Q4}
\end{align}

\subsection{Calculation of coefficients in Oseen-Frank energy deduced from the tensor model}

We denote
\begin{align*}
I_1=(\nabla\cdot\nn)^2,\quad I_2=(\nn\cdot\nabla\times\nn)^2,\quad I_3=|\nn\times(\nabla\times\nn)|^2,
\quad I_4=\mathrm{tr}(\nabla\nn)^2-(\nabla\cdot\nn)^2.
\end{align*}
Then we can easily get that:
\begin{align}
(\partial_in_i)^2=I_1,\quad (n_i\partial_in_k)^2=I_3,\quad \partial_in_j\partial_jn_i=I_1+I_4,\quad (\partial_in_j)^2=I_1+I_2+I_3+I_4.
\end{align}
From the identity
\begin{align}
Q_{4jklm}=&S_4\Big(n_in_kn_jn_l-\frac17\big(n_in_j\delta_{kl}+n_in_k\delta{jl}
+n_in_l\delta_{jk}\nonumber\\
&+n_kn_l\delta_{ij}+n_jn_l\delta_{ik}+n_jn_k\delta_{il}\big)
+\frac{1}{35}\big(\delta_{ij}\delta_{kl}+\delta_{il}\delta_{jk}+\delta_{ik}\delta_{jl}\big)\Big),
\end{align}
and the fact that $Q_{ijkk}=Q_{ijjk}=\cdots=0$, we have
\begin{align}
|\nabla(cQ_4)|^2=&S_4^2\partial_h\Big(n_in_kn_jn_l-\frac17\big(n_in_j\delta_{kl}+n_in_k\delta_{jl}
+n_in_l\delta_{jk}+n_kn_l\delta_{ij}\nonumber\\
&+n_jn_l\delta_{ik}+n_jn_k\delta_{il}\big)
+\frac{1}{35}\big(\delta_{ij}\delta_{kl}+\delta_{il}\delta_{jk}
+\delta_{ik}\delta_{jl}\big)\Big)\partial_h(n_in_kn_jn_l)\nonumber\\
=&4S_4^2\partial_h\Big(n_in_kn_jn_l-\frac17\big(n_in_j\delta_{kl}+n_in_k\delta_{jl}
+n_in_l\delta_{jk}+n_kn_l\delta_{ij}\nonumber\\
&+n_jn_l\delta_{ik}+n_jn_k\delta_{il}\big)
+\frac{1}{35}\big(\delta_{ij}\delta_{kl}+\delta_{il}\delta_{jk}
+\delta_{ik}\delta_{jl}\big)\Big)n_kn_jn_l\partial_hn_i\nonumber\\
=&4S_4^2\Big((\partial_hn_i)^2-\frac37(\partial_hn_i)^2\Big)=\frac{16}{7}S_4^2(I_1+I_2+I_3+I_4).
\end{align}
Similarly, we can obtain
\begin{align*}
&\partial_i(cQ_{4iklm})\partial_j(cQ_{4jklm})=
S_4^2\big(\frac{46}{49}I_1+\frac{6}{49}I_2+\frac{60}{49}I_3+\frac{12}{49}I_4\big),\\
&\partial_i(cQ_{4jklm})\partial_j(cQ_{4iklm})=
S_4^2\big(\frac{46}{49}I_1+\frac{6}{49}I_2+\frac{60}{49}I_3+\frac{40}{49}I_4\big),\\
&\partial_i(cQ_{4ijkl})\partial_j(cQ_{2kl})=
S_4S_2\Big(\frac{8}{7}I_3-\frac{2}{7}(3I_1+I_2+2I_4)\Big),\\
&|\nabla(cQ_2)|^2=2S_2^2(I_1+I_2+I_3+I_4),\\
&\partial_i(cQ_{ik})\partial_j(cQ_{jk})=S_2^2(I_1+I_3),\\
&\partial_i(cQ_{jk})\partial_j(cQ_{ik})=S_2^2\Big(I_1+I_3+I_4\Big).
\end{align*}

Substituting the above equalities into (\ref{Q-Energy}), we can get that
\begin{align}
F_{\text{elastic}}^{(2)}
=&\frac{c^2}2\int_\Omega\bigg\{2J_2S_2^2(I_1+I_2+I_3+I_4)
+J_3\frac{16}{7}S_4^2(I_1+I_2+I_3+I_4)\nonumber\\
&+J_5\Big(S_2^2(I_1+I_3)+S_2^2(I_1+I_3+I_4)\Big)+J_7S_2S_4\Big(\frac{8}{7}I_3-\frac{2}{7}(3I_1+I_2+2I_4)\Big)\nonumber\\
&+J_6\Big(S_4^2\big(\frac{46}{49}I_1+\frac{6}{49}I_2+\frac{60}{49}I_3+\frac{12}{49}I_4\big)
+S_4^2\big(\frac{46}{49}I_1+\frac{6}{49}I_2+\frac{60}{49}I_3+\frac{40}{49}I_4\big)
\Big)\bigg\}\ud\xx.\nonumber
\end{align}
Therefore, we have
\begin{align*}
&K_1=c^2\Big(2S_2^2(J_2+J_5)+S_4^2(\frac{16}{7}J_3+\frac{92}{49}J_6)-\frac67J_7S_2S_4\Big),\\
&K_2=c^2\Big(2S_2^2J_2+S_4^2(\frac{16}{7}J_3+\frac{12}{49}J_6)-\frac27J_7S_2S_4\Big),\\
&K_3=c^2\Big(2S_2^2(J_2+J_5)+S_4^2(\frac{16}{7}J_3+\frac{120}{49}J_6)+\frac87J_7S_2S_4\Big),\\
&K_4=c^2\Big(S_2^2J_5+\frac{40}{49}J_6S_4^2-\frac27J_7S_2S_4\Big).
\end{align*}
Substituting the expressions of $J_i$:
{\begin{align*}
&J_2=-\frac{\pi}{2}{L^5\eta}k_BT\Big(\frac32\alpha_{12}+\frac37\alpha_{22}
+\frac{9}{49}\alpha_{32}\Big),\qquad
J_3=-\frac{35\pi}{16}{L^5\eta}k_BT\alpha_{13},\qquad\\
&J_5=-\frac{\pi}{2}{L^5\eta}k_BT\Big(\alpha_{31}+\frac{25}{49}\alpha_{32}+\frac67\alpha_{22}\Big),\\
&J_6=-\frac{3\pi}{4}{L^5\eta}k_BT\alpha_{32},\qquad
J_7=-\frac{\pi}{2}{L^5\eta}k_BT\big(3\alpha_{22}+\frac{18}{7}\alpha_{32}\big).
\end{align*}}
into it, the elastic coefficients can be written as
\begin{align*}
&K_1=\pi c^2L^5\eta k_BT\Big(-S_2^2\Big(\eta^2\frac{299}{1568}-\frac{15}{7\cdot64}
-\frac{12\ln2\eta^2}{49}\Big)+S_4^2\Big(\frac{15\eta^2}{128}
-\frac{115\eta^2}{49}\Big(\frac{1}{8}-\frac{\ln2}{2}\Big)\Big)\nonumber\\
&\qquad\qquad+S_2S_4\frac{15}7\Big(-\frac{1}{64}+\frac{5\eta^2}{14}-\frac{3}{7}\eta^2{\ln2}\Big)\Big),\\
&K_2=\pi c^2{L^5\eta}k_BT\Big(-5S_2^2\Big(\eta^2\frac{19}{1568}-\frac{1}{7\cdot64}
-\frac{3\ln2\eta^2}{98}\Big)+S_4^2\Big(\frac{15\eta^2}{128}-\frac{15}{49}\eta^2\Big(\frac{1}{8}-\frac{\ln2}{2}\Big)\Big)\\
&\qquad\qquad+S_2S_4\frac{5}{7}\Big(-\frac{1}{64}+\frac{5\eta^2}{14}-\frac{3}{7}\eta^2{\ln2}\Big)\Big),\\
&K_3=\pi c^2{L^5\eta}k_BT\Big(-S_2^2\Big(\eta^2\frac{299}{1568}-\frac{15}{7\cdot64}
-\frac{12\ln2\eta^2}{49}\Big)+S_4^2\Big(\frac{15\eta^2}{128}
-\frac{150}{49}\eta^2\Big(\frac{1}{8}-\frac{\ln2}{2}\Big)\Big)\\
&\qquad\qquad-S_2S_4\frac{20}7\Big(-\frac{1}{64}+\frac{5\eta^2}{14}-\frac{3{\ln2}\eta^2}{7}\Big)\Big),\\
&K_4=\pi c^2{L^5\eta}k_BT\Big(\frac{S_2^2}{2}\Big(\frac{9\ln2\eta^2}{98}-\frac{51\eta^2}{392}+\frac{5}{7\cdot32}\Big)
-S_4^2\frac{25\eta^2}{49}\Big(\frac{1}{4}-{\ln2}\Big)\\
&\qquad\qquad+S_2S_4\frac{5}7\Big(-\frac{1}{64}+\frac{5\eta^2}{14}-\frac{3{\ln2}\eta^2}{7}\Big)\Big).
\end{align*}

\bigskip

\noindent {\bf Acknowledgments.} The authors are  grateful to Prof. Zhifei Zhang
for his suggestions which helped improve the paper greatly. The authors also would like to thank
Dr. Weiquan Xu for his help on numerical simulation.
P. Zhang is partly supported by NSF of China under Grant 50930003 and 21274005.

\end{document}